\documentclass[number,citesort,MSNbibl,dvips]{arxbj}
\usepackage{upgreek}

\aid{0}
\volume{18}
\issue{2}
\pubyear{2012}
\firstpage{703}
\lastpage{734}
\doi{10.3150/10-BEJ347}

\makeatletter
\newcommand{\E}{\mathbb{E}}
\renewcommand{\P}{\mathbb{P}}
\newcommand{\N}{\mathbb{N}}
\newcommand{\R}{\mathbb{R}}
\newcommand{\C}{\mathbb{C}}
\newcommand{\cB}{\mathcal{B}}
\newcommand{\cD}{\mathcal{D}}
\newcommand{\cE}{\mathcal{E}}
\newcommand{\cF}{\mathcal{F}}
\newcommand{\cL}{\mathcal{L}}

\newcommand{\cR}{\mathcal{R}}
\newtheorem{theo}{Theorem}
\newtheorem{pro}{Proposition}
\newtheorem{lem}{Lemma}
\newtheorem{alem}{Lemma}[section]
\newremark{rem}{Remark}
\newproclaim{assumptions}{Assumptions}
\newproclaim{assumption}{Assumption}
\newproclaim{cnd}{Condition (${D_{m_0}}$)}
\newcommand{\cal}{\mathcal}
\newcommand{\fraca}[2]{{#1}/{#2}}
\newcommand{\fracc}[2]{{#1}/{(#2)}}
\newcommand{\fracb}[2]{{(#1)}/{#2}}
\makeatother

\begin{document}
\begin{frontmatter}

\title{A uniform Berry--Esseen theorem on $M$-estimators for
geometrically ergodic~Markov chains}
\runtitle{Berry--Esseen theorem on $M$-estimators for geometrically
ergodic Markov chains}

\begin{aug}
\author[1]{\fnms{Lo\"{\i}c} \snm{Herv\'e}\thanksref{1,e1}\ead[label=e1,mark]{loic.herve@insa-rennes.fr}},
\author[1]{\fnms{James} \snm{Ledoux}\thanksref{1,e2}\ead[label=e2,mark]{james.ledoux@insa-rennes.fr}} \and
\author[1,2]{\fnms{Valentin} \snm{Patilea}\corref{}\thanksref{1,2}\ead[label=e3]{valentin.patilea@insa-rennes.fr}}
\runauthor{L. Herv\'e, J. Ledoux and V. Patilea}
\address[1]{Universit\'e Europ\'eenne de Bretagne, INSA-IRMAR, UMR-CNRS
6625. Institut National des Sciences Appliqu\'{e}es de Rennes, 20,
Avenue des Buttes de Co\"{e}smes CS 14315,
35043 Rennes cedex, France. \printead{e1,e2}}
\address[2]{CREST (Ensai), Campus de Ker-Lann,
Rue Blaise Pascal -- BP 37203, 35172 Bruz cedex, France. \printead{e3}}
\end{aug}

\received{\smonth{9} \syear{2009}}
\revised{\smonth{9} \syear{2010}}

%
\begin{abstract}
Let $\{X_n\}_{n\ge0}$ be a $V$-geometrically ergodic Markov chain.
Given some real-valued functional~$F$, define $M_n(\alpha) := n^{-1}
\sum_{k=1}^n
F(\alpha,X_{k-1},X_k)$, $\alpha\in\mathcal{A}\subset\R$. Consider
an $M$ estimator
$\widehat\alpha_n $, that is, a~measurable function of the
observations satisfying $M_n(\widehat\alpha_n) \leq
\min_{\alpha\in\mathcal{A}} M_n(\alpha) + c_n$ with $\{c_n\}
_{n\geq1}$
some sequence of real numbers going to zero. Under some
standard regularity and moment assumptions, close to those of the i.i.d. case,
the estimator $\widehat\alpha_n $ satisfies a Berry--Esseen theorem
uniformly with respect to the underlying probability distribution of
the Markov chain.
\end{abstract}

%

\vspace*{-3pt}
\end{frontmatter}

\section{Introduction}\label{intro}

Let $(E,\cE)$ be a measurable space with $\cE$ a countably generated
$\sigma$-field, and let $\{X_n\}_{n\ge0}$ be a Markov chain
with state space $E$ and transition kernels $\{Q_{\theta}(x,\cdot)
\dvt
x\in E\},$ where~$\theta$ is a parameter in some general set
$\Theta$. The initial distribution of the chain, that is, the
probability distribution of~$X_0$, is denoted by~$\mu$ and
may or may not depend on~$\theta$.
Although $\{X_n\}_{n\geq0}$
does not need to be the canonical version, we use the standard
notation~$\P_{\theta,\mu}$ to refer to the probability distribution
of $\{X_n\}_{n\geq0}$ (and
$\E_{\theta,\mu}$ for the expectation w.r.t.~$\P_{\theta,\mu}$).
We consider that $\{X_n\}_{n\geq0}$ is a $V$-geometrically ergodic
Markov chain,
where $V \dvtx E{\rightarrow}[1,+\infty)$ is some fixed \textit
{unbounded} function.
This class of Markov chains is large enough to cover interesting
applications (see \cite{MeyTwe93}, Sections 16.4 and 16.5).

The parameter of interest is $\alpha_0 = \alpha_0(\theta)\subset
\mathcal{A}$, where $\alpha_0(\cdot)$ is a function of the
parameter~$\theta$ and $\mathcal{A}$ is an open interval of $\R$. To estimate
$\alpha_0$, let us
introduce the statistic
%
\begin{equation}\label{def_m}
M_n(\alpha) := \frac{1}{n} \sum_{k=1}^n F(\alpha,X_{k-1},X_k),\vadjust{\goodbreak}
\end{equation}
where $F$ is a real-valued measurable functional
on $\mathcal{A}\times E^2 $. We define an \textit{$M$-estimator}
(this is slightly more general than the usual definition of
$M$-estimators or minimum contrast estimators, where $c_n = 0$, see
\cite{arc98}) to be a
random variable $\widehat\alpha_n$ depending on the observations
$(X_0,\ldots,X_n)$
such that
\[
M_n(\widehat\alpha_n )\leq\min_{\alpha\in\mathcal{A}} M_n(\alpha)
+c_n,
\]
where $\{c_n\}_{n\geq1}$ is a sequence of non-negative real numbers
going to zero to be specified later. Assume that for all $\theta\in
\Theta$
\[
M_{\theta} (\alpha):= \lim_{n\rightarrow\infty} \E_{\theta,\mu
}[M_n (\alpha)]
\]
is well defined everywhere on $\mathcal{A}$ and does not depend on
$\mu$. In addition, assume that there exists a unique ``true'' value
$\alpha_0$ of
the parameter of interest, that is, $M_{\theta}(\alpha_0) <
M_{\theta}(\alpha)$, $\forall\alpha\neq\alpha_0$.
We want to prove the following \textit{uniform} Berry--Esseen bound for
$\widehat\alpha_n$
%
\renewcommand{\theequation}{BE}
\begin{equation} \label{BE}
\sup_{\theta\in\Theta}\sup_{u\in\R} \biggl|\P_{\theta,\mu
} \biggl\{\frac{\sqrt n}{\tau(\theta)}  (\widehat{\alpha}_n-\alpha
_0) \leq u \biggr\}
- \Gamma(u) \biggr| = \mathrm{O}\biggl(\frac{1}{\sqrt n} \biggr),
\end{equation}%
\renewcommand{\theequation}{\arabic{equation}}\setcounter{equation}{1}%
where $\Gamma$ denotes the standard normal distribution function,
and $\tau(\theta)$ is some positive real number defined in Theorem
\ref{Th-B-E-estimator}.

To derive (\ref{BE}), we use Pfanzagl's approach \cite{Pfa71}.
Besides technical assumptions, this approach relies on
several ingredients. First, we need the uniform
consistency condition:
\begin{enumerate}[(UC)]
\item[(UC)] $\forall d>0$, $\sup_{\theta\in\Theta} \P
_{\theta,\mu}  \{| \widehat\alpha_n -\alpha_0| \geq d  \}
= \mathrm{O}(1/\sqrt{n} )$.
\end{enumerate}
Second, consider the following two convergence properties: If
$S_n(\alpha_0) := \sum_{k=1}^{n} \xi(\alpha_0,\break X_{k-1},X_k)$ with
$\xi(\alpha_0, X_{k-1},X_k)$ centered,
\begin{enumerate}[(b)]
\item[(a)]\hypertarget{a}{}\textit{the sequence $ \{\E_{\theta,\mu}
[S_n^2(\alpha_0)]/n \}_{n\ge1}$ converges to a real number
$\sigma^2(\theta)$};
\item[(b)]\hypertarget{b}{} \textit{there exists a positive constant $B(\xi)$ such
that for any $n\geq1$}
\[
\sup_{\theta\in\Theta}\sup_{u\in\R} \biggl|\P_{\theta,\mu
} \biggl\{\frac{S_n(\alpha_0)}{\sigma(\theta)  \sqrt n}
\leq u \biggr\}- \Gamma(u) \biggr| \leq\frac{B(\xi)}{\sqrt n}.
\]
\end{enumerate}
The properties ({a}) and ({b}) will be required for certain
$\xi(\alpha_0, x,y)$ defined as
linear combinations of some functionals related to $F$.
To obtain (a) and (b) for such $\xi(\alpha_0, x,y)$ with
$V$-geometrically ergodic Markov chains, a natural moment (or
$V$-domination) condition is used: There exist positive
constants $C_\xi$ and $m$ such that
%
\begin{equation} \label{M}
\forall(x,y)\in E^2,  \forall\alpha\in\mathcal{A} \qquad
|\xi(\alpha, x,y) |^m \leq C_\xi  \bigl(V(x)+V(y)\bigr).
\end{equation}

The paper is organized as follows. In Section \ref{IV.1}, an extended
version of Pfanzagl's theorem \cite{Pfa71},  is stated for
any sequence of observations, not necessarily Markovian. Section~\ref
{IIbis} is devoted to a Berry--Esseen bound for the additive functional
$\sum_{k=1}^n\xi(\alpha_0,X_{k-1},X_k)$ of a $V$-geometrically
ergodic Markov chain $\{X_n\}_{n\ge0}$ with $\xi$ satisfying
inequality (\ref{M}).
In Section~\ref{II}, we prove that the properties (a) and (b) are
fulfilled when inequality~(\ref{M}) holds with the (almost expected)
order $m$, namely: $m>2 \Rightarrow$ ({a}), and $m>3
\Rightarrow$ ({b}). These results\vadjust{\goodbreak} follow from the weak spectral
method based on the theorem of Keller and Liverani~\cite{KelLiv99}.
This approach, introduced in \cite{HenHer04}, is fully described in
\cite{HerPen08} in the Markov context (see also \cite{GouLiv06,Gou08}
and other references given in \cite{HerPen08}).
It~is important to notice that Pfanzagl's method requires the precise
control of the constant $B(\xi)$ in property ({b}) as a~function
of the size of $\xi$. The present operator-type approach shows that
$B(\xi)$ depends only on the constant $C_\xi$ in inequality (\ref
{M}). Thanks to these preliminary results, in Section~\ref{III.2}
we prove our main statement, that is:
\begin{enumerate}[(R)]
\item[(R)] \textit{Under some technical assumptions and the uniform
consistency condition \textup{(UC)}, if two functionals $F'$ and $F''$
related to $F$ (in the basic case $F'$ and $F''$ are the first- and
second-order derivatives of $F$ with respect to $\alpha$) satisfy
inequality  \textup{(\ref{M})} for some $m>3$ and constants $C_{F'}$,
$C_{F''}$ that do not depend on $\alpha$, then $\widehat{\alpha}_n$
satisfies property \textup{(\ref{BE})}.}
\end{enumerate}
To the best of our knowledge, the result  {(R)} is new. It
completes the central limit theorem for $\{\widehat{\alpha}_n\}_{n\ge
1}$ proved in \cite{DehYao07} when inequality (\ref{M}) holds with
$m=2$. The domination condition~(\ref{M}) required by  {(R)} is
almost optimal in the sense that we impose $m>3$ in place of the best
possible value $m=3$ obtained in the i.i.d. case. In Section \ref
{iid}, our results are applied to the AR(1) process with ARCH
(autoregressive conditional heteroscedastic) of order-1 errors. The
paper ends with a conclusion section.

Let us close the \hyperref[intro]{Introduction} with a brief review of previous related
works in the literature.
In \cite{Pfa71}, $\{X_n\}_{n\in\N}$ is a sequence of i.i.d. random
variables and Pfanzagl proved a Berry--Esseen theorem for minimum
contrast estimators (which are special instances of $M$-estimators)
associated with functionals of the form $F(\alpha,X_k)$.
In \cite{Pfa71}, the moment conditions on $F':=\partial F/\partial
\alpha$, $F'':=\partial^2 F/\partial\alpha^2$ are the expected ones
since the property (b) is fulfilled under the expected third moment
condition \cite{Fel71}, Chapter~XVI. Using convexity arguments,
Bentkus \textit{et al.} \cite{BenBloGot97} proposed an alternative
method for deriving Berry--Esseen bounds for $M$-estimators with i.i.d. data.
In the Markov context, the method proposed by Pfanzagl is extended,
first by Rao to cover the case of uniformly ergodic Markov chains \cite
{Rao73}, second in \cite{MilRau89} to the case of the linear
autoregressive model. However, their assumptions to get the property (BE)
include much stronger moment conditions involving both the functional
$F$ and the Markov chain.
Here, as already mentioned, the weak spectral method of \cite
{HerPen08} enables us to have an (almost) optimal treatment of (a) and
(b), and hence an improved Berry--Esseen result (BE).

\section{The Pfanzagl method revisited} \label{IV.1}

We state and prove a general result that allows us to derive uniform
Berry--Esseen
bounds for $M$-estimators. This result is an extended version of
Theorem 1 in \cite{Pfa71} and is applied to our Markov context in
Section~\ref{III.2}.

\subsection{The result}
Consider a statistical model $ (\Omega, \cF,\{\P_{\theta},
\theta\in
\Theta\} )$, where $\Theta$ denotes some parameter space, and let
$\{X_n\}_{n\geq0}$ be any sequence of observations (not necessarily
Markovian). Let us denote the expectation with respect to
$\P_{\theta}$ by $\E_{\theta}$.\vadjust{\goodbreak}

For each $n$, let $M_n(\alpha) $ be a measurable functional of the
observations $X_0,\ldots, X_n$ and the parameter of interest $\alpha
\in\mathcal{A}$, where $\mathcal{A}$ is some open
interval of $\R$.
Let $\{c_n\}_{n\ge1}$ be a~sequence of non-negative real numbers going to
zero at some rate to be specified later. An \textit{$M$-estimator}
is a measurable function $\widehat\alpha_n$ of the observations
$(X_0,\ldots,X_n)$ such that
%
\begin{equation} \label{Def_Mn}
M_n(\widehat\alpha_n)\leq \min_{\alpha\in\mathcal{A}} M_n(\alpha)
+ c_n.
\end{equation}
This is the usual definition of
minimum contrast estimators as soon as $c_n\equiv0$.\vspace*{-3pt}
\begin{assumptions*}
 Suppose that for all $n\geq1$ and
$\alpha\in\mathcal{A}$,
there exist $M'_n(\alpha)$,
$M''_n(\alpha)$ some measurable functions depending on $X_0, X_1,\ldots,
X_n$ and on the parameter of interest, such that the following
properties hold true:
\begin{enumerate}[(A4$'$)]
\item[(A1)] \label{conda1} $\forall\theta\in\Theta$, there exists
a unique $\alpha_0 = \alpha_0 (\theta)\in
\mathcal{A}$ such that $ M'_\theta(\alpha_0) = 0$, where $ M'_\theta
(\alpha) := \lim_{n\rightarrow\infty} \E_{\theta}[ M'_n(\alpha
)]$ (the limit is
assumed to be well defined for all $(\theta,\alpha)\in\Theta\times
\mathcal{A}$);
\item[(A2)] $0 < \inf_{\theta\in\Theta} m(\theta)
\leq\sup_{\theta\in\Theta}m(\theta) < \infty$, where $m(\theta)
:= \lim_{n\rightarrow\infty} \E_{\theta}[ M''_n(\alpha_0)]$ (the
limit is assumed to be well defined for all $\theta$);
\item[(A3)] for every $n\geq1$, there exists $r_n >0$ independent of
$\theta$
such that $r_n=\mathrm{o}(n^{-1/2})$~and
\[
\sup_{\theta\in\Theta} \P_\theta  \{|M'_n(\widehat{\alpha
}_n)|\geq r_n \} = \mathrm{O}(n^{-1/2});
\]
\item[(A4)] for $j=1,2$, there exists a function $\sigma_j(\cdot)$
such that
$0< \inf_{\theta\in\Theta} \sigma_j(\theta) \leq\break
\sup_{\theta\in\Theta}\sigma_j(\theta)< \infty$ and there exists
a positive constant $B$ such that for all $n\geq1$
\begin{eqnarray*}
\sup_{\theta\in\Theta}\sup_{u\in\R} \biggl|\P_{\theta} \biggl\{
\frac{\sqrt n}{\sigma_1(\theta)}  M'_n(\alpha_0) \leq u \biggr\} -
\Gamma(u) \biggr|
&\leq&\frac{B}{\sqrt n},
\\[-2pt]
\sup_{\theta\in\Theta}\sup_{u\in\R} \biggl|\P_{\theta}  \biggl\{
\frac{\sqrt n}{\sigma_2(\theta)} \bigl(M''_n(\alpha_0)-m(\theta)
 \bigr)
\leq u \biggr\} - \Gamma(u) \biggr|& \leq&\frac{B}{\sqrt n};
\end{eqnarray*}
\item[(A4$'$)] \textit{for $n\geq1$, $|u| \leq2\sqrt{\ln n}$ and
$\theta\in\Theta$, there is a positive number $\sigma_{n,u}(\theta
)$ such that}
\begin{eqnarray*}
|\sigma_{n,u}(\theta) - \sigma_1(\theta)| &\leq& A'  \frac
{|u|}{\sqrt n},\\[-2pt]
 \biggl|\P_{\theta} \biggl\{\frac{\sqrt{n}}{\sigma_{n,u}(\theta
)} \biggl(M'_n(\alpha_0) +
\frac{u  \sigma_1(\theta)}{\sqrt{n} m(\theta)
} \bigl(M''_n(\alpha_0)-m(\theta)  \bigr) \biggr) \leq u \biggr\} -
\Gamma(u) \biggr| &\leq&\frac{B'}{\sqrt n}
\end{eqnarray*}
with some positive constants $A',B'$ independent of $n,u,\theta$;
\item[(A5)] for any $(\alpha, \tilde\alpha)\in\mathcal{A}^2$, let
$R_n(\alpha, \tilde\alpha)$ be defined by the equation
\[
M'_n(\tilde\alpha) = M'_n(\alpha) + [M''_n(\alpha)
+ R_n(\alpha, \tilde\alpha)](\tilde\alpha- \alpha).
\]
For each $n$, there exist $\omega_n \geq 0$ and a real-valued
measurable function $W_n $ depending on $X_0,\ldots,X_n$, both
independent of $\theta$, such that $\omega_n =\mathrm{o}(1)$ and
\[
\forall
(\alpha,\tilde\alpha)\in\mathcal{A}^2 \qquad   | R_n(\alpha,
\tilde\alpha)|\leq \{|\alpha-\tilde\alpha|+ \omega_n \}
  W_n,\vadjust{\goodbreak}
\]
and there is a constant $c_W > 0$ such that
\[
\sup_{\theta\in\Theta} \P_\theta\{ c_W \leq W_n \} = \mathrm{O}(n^{-1/2});
\]
\item[(A6)] $\widehat{\alpha}_n$ is assumed to be uniformly
consistent, that is, there exists $\gamma_n=\mathrm{o}(1)$ such that
\[
\sup_{\theta\in\Theta}\P_{\theta} \{
|\widehat{\alpha}_n-\alpha_0| \geq d  \} \leq\gamma_n ,
\]
where $d:=\inf_{\theta\in\Theta} m(\theta) /8c_W$ with $c_W$ and
$m(\theta)$ defined
in  \textup{(A5)} and \textup{(A2)}, respectively.
\end{enumerate}
\end{assumptions*}

Let us comment on these assumptions. Condition (A1)
identifies the true value of the parameter. In conditions (A1) and
(A2), the expectations $\E_\theta[M'_n(\alpha)]$ and
$\E_\theta[M''_n(\alpha_0)]$ may depend on $n$, as in the Markovian
framework considered in the sequel when the
initial distribution is not the stationary distribution. Condition
(A3) ensures that the estimator (approximately) satisfies a kind of
first-order condition. Such a condition allows us to
take into account the numerical errors with which we are faced when
computing~$\widehat\alpha_n$. It may also be useful when the estimator of the
parameter $\alpha_0$ depends on some
``nuisance'' parameters (see the example in the second part of
Section~\ref{iid}).
Conditions (A4) and (A4$'$) are the uniform Berry--Esseen bounds for
$M'_n(\alpha_0)$, $M''_n(\alpha_0)$ and for some of their linear combinations.
The identity defining $R_n(\alpha,\tilde\alpha)$ in condition (A5) is
guaranteed by a Taylor expansion when the criterion $M_n(\alpha)$ is
twice differentiable with respect to~$\alpha$. In this case~$M'_n$
and $M''_n$ are nothing else but the first- and second-order
derivatives of~$M_n$ with respect to $\alpha$.
The reminder $R_n(\alpha,\tilde\alpha)$ must satisfy a Lipschitz
condition. For instance, when
$\omega_n = 0$, this holds true if $\alpha\mapsto M_n(\alpha)$ is
three times continuously differentiable with a bounded third-order
derivative. Condition (A6) is a standard consistency condition
(see~\cite{BenBloGot97}). General sufficient conditions for (A6) with
$\gamma_n = \mathrm{O}(n^{-1})$ have been proposed in the case of
i.i.d. observations or uniformly ergodic Markov chains (see \cite
{MicPfa71}, Lemma 4, and \cite{Rao73}, Lemma 4.1, resp.). Such
general arguments can easily be adapted to the geometrically ergodic
Markov chain framework. In specific examples, like the one investigated
in Section~\ref{iid}, condition (A6) can be checked by direct arguments.

The proof of Theorem~\ref{th-pfanzagl}, which adapts the arguments of
\cite{Pfa71}, is given in Section~\ref{ProofTh1}.

\begin{theo} \label{th-pfanzagl} Under conditions  \textup{(A1)--(A6)},
there exists a positive constant $C$ such that
%
\begin{equation} \label{Bth1}
\forall n\geq1 \qquad    \sup_{\theta\in\Theta}\sup_{u\in\R}
 \biggl|\P_{\theta} \biggl\{\frac{\sqrt n}{\tau(\theta)}
(\widehat{\alpha}_n-\alpha_0) \leq u \biggr\} - \Gamma(u) \biggr|
\leq
C \biggl( \frac{1}{\sqrt n} + \sqrt{n} r_n + \omega_n +
\gamma_n \biggr)
\end{equation}
with $\tau(\theta) :=
\sigma_1(\theta)/m(\theta) $.
\end{theo}

To obtain the classical order $\mathrm{O}(n^{-1/2})$ of the Berry--Esseen bound,
one needs $\gamma_n=\mathrm{O}(n^{-1/2})$, $r_n=\mathrm{O}(n^{-1})$ and $\omega_n
=\mathrm{O}(n^{-1/2})$. Note that this usually requires that the~se\-quence $\{
c_n\}_{n\ge1}$ in (\ref{Def_Mn}) decreases at the rate $n^{-3/2}$.
This is to be compared to the ra\-te~$n^{-1}$ that is usually required to
obtain the asymptotic normality of $M$-estimators (see \cite{arc98}).
\begin{rem} A close inspection of the proof of Theorem~\ref
{th-pfanzagl} below shows that the constant~$C$ in inequality (\ref
{Bth1}) can be tracked provided that the $\mathrm{O}(\cdot)$ and $\mathrm{o}(\cdot)$
rates in assumptions~(A3)--(A6)  are more explicit. For the sake
of brevity, we only consider the case where $c_n=r_n=\omega_n=0$,
$\alpha(\theta)=\theta$ and  (A3)  is: for any $n\ge1$,
$|M'_n(\widehat{\theta}_n)| = 0$. The constants~$C$ in the various
inequalities of assumptions  (A4)--(A6)  are denoted by $C_1,C_2$
in  {(A4)}, $C_3,C_4$ in  (A4$'$) and $C_5$ in  {(A5)}
and we choose $\gamma_n \leq C_6   n^{-1/2}$ in  {(A6)}. Then
we can obtain from Propositions \ref{Prop1} and \ref{Prop2} that
\[
\label{Ineg_BE_constante}
\forall n \ge1 \qquad   \sup_{\theta\in\Theta} \biggl|\P_{\theta
} \biggl\{\frac{\sqrt
n}{\tau(\theta)}  (\widehat{\alpha}_n - \alpha_0) \leq u \biggr\}
- \Gamma(u) \biggr| \leq\frac{C}{\sqrt n},
\]
 where  $C:= \frac{1}{2} + \frac{1}{\sqrt{2\uppi
}} + 2C_1 + 2C_2 + \frac{\exp(-a^2/2)}{a} + C_5 + C_6$   when
$|u| \geq2\sqrt{\ln n}$;  or
  $C:= 2 [\frac{1}{\sqrt{2\uppi}} + 2C_1 + 4C_2 +
2\frac{\exp(-a^2/2)}{a} + 2C_5 + C_6 ] + C_4+ \frac{16
\mathrm{e}^{-1}(C_3+\overline{\sigma}^2 c_W)}{\underline{\sigma}_1 \sqrt
{2\uppi}}$   when   $|u| < 2\sqrt{\ln n}$   provided that   $\sqrt
{n/\ln n} \ge\max (8 c_W \overline{\sigma}^2,4
)/\underline{\sigma}_1;$
with $a:=\inf_{\theta\in\Theta} (m(\theta)/4\sigma_2(\theta
) )$, $\overline{\sigma}:=\sup_{\theta\in\Theta} \sigma
_1(\theta)/m(\theta)$, $ \underline{\sigma}_1 :=
\inf_{\theta\in\Theta}\sigma_1(\theta)$.
\end{rem}

\subsection{\texorpdfstring{Proof of Theorem~\protect\ref{th-pfanzagl}}{Proof of Theorem 1}} \label{ProofTh1}

The hypotheses of
Theorem~\ref{th-pfanzagl} are assumed to hold. For the sake of
brevity, the sequence $\{r_n\}_{n\ge1}$ in (A3) is supposed to be
such that $r_n=\mathrm{o}(n^{-1/2})$
and $|M'_n(\widehat\alpha_n)| \leq r_n$ for every $n\geq1$. In the
general case, it suffices to work on the event
$\{|M'_n(\widehat\alpha_n)| \leq r_n\}$ and to bound the probability
of the event $\{|M'_n(\widehat\alpha_n)| > r_n\}$ using (A3). From
conditions (A2) and (A4),
\begin{eqnarray*}
\tau(\theta) &:=& \frac{\sigma_1(\theta)}{m(\theta)}, \qquad   \underline
{m} :=\inf_{\theta\in\Theta} m(\theta), \qquad   \overline{m} := \sup
_{\theta\in\Theta}m(\theta), \\
 \underline{\sigma}_j &:=& \inf_{\theta\in\Theta} \sigma_j(\theta
), \qquad
 \overline{\sigma}_j := \sup_{\theta\in\Theta}\sigma_j(\theta),
\end{eqnarray*}
$j=1,2,$ are well defined.
Recall that $0 <\underline{m} \leq\overline{m} < \infty$ and
$0<\underline{\sigma}_j \le\overline{\sigma}_j< \infty$.
Note that the function $\tau(\cdot)$ is positive and bounded.
In the following, $C$ denotes a positive constant whose value may be
different from line to line.

Inequality (\ref{Bth1}) is proved, first for $|u|\geq2\sqrt{\ln n}$,
second for $|u|< 2\sqrt{\ln n}$. In fact, for $|u|\geq2\sqrt{\ln
n}$, the bound in inequality (\ref{Bth1}) does not involve $r_n$ and
$\omega_n$.
\begin{pro} \label{Prop1} There exists a positive constant $C$ such
that for each
$n\geq1$ and all $u\in\R$ such that $|u| \geq2\sqrt{\ln n}$
%
\begin{equation} \label{Ineg_BE}
\sup_{\theta\in\Theta} \biggl|\P_{\theta} \biggl\{\frac{\sqrt
n}{\tau(\theta)}  (\widehat{\alpha}_n - \alpha_0) \leq u \biggr\}
- \Gamma(u) \biggr| \leq\frac{C}{\sqrt n} + \gamma_n.
\end{equation}
\end{pro}

\begin{pf} For $|u|\ge 2\sqrt{\ln n}$, it is easily checked that
\[
 \biggl|\P_{\theta} \biggl\{\frac{\sqrt n}{\tau(\theta)}  (\widehat
{\alpha}_n - \alpha_0) \leq u \biggr\}
- \Gamma(u) \biggr|
\le \P_{\theta} \biggl\{\frac{\sqrt n}{\tau(\theta)}   |
\widehat{\alpha}_n - \alpha_0  | \geq2\sqrt{\ln n}  \biggr\} +
\Gamma\bigl(-2\sqrt{\ln n}\bigr).
\]
Now,
\[
\Gamma\bigl(-2\sqrt{\ln n}\bigr) \leq \frac{1}{2\sqrt{\ln n}}   \frac
{1}{\sqrt{2\uppi}} \int_{2\sqrt{\ln n}}^{+\infty} v
\mathrm{e}^{-\fraca
{v^2}{2}}\,\mathrm{d}v =
\frac{1}{2\sqrt{\ln n}}   \frac{1}{\sqrt{2\uppi}} \frac{1}{n^2}.
\]
Finally, the proof is complete if there exists $C>0$ such that (see
\cite{MicPfa71}, Lemma 6)
%
\begin{equation} \label{Lemme6}
\forall n\geq1 \qquad   \sup_{\theta\in\Theta}\P_{\theta}\biggl \{
\frac{\sqrt
n}{\tau(\theta)}   | \widehat{\alpha}_n-\alpha_0  | >
2\sqrt{\ln n}  \biggr\} \leq\frac{C}{\sqrt n} + \gamma_n.
\end{equation}

It follows from (A5) and (A3) that $ |M'_n(\alpha_0)| + r_n \geq
|\widehat{\alpha}_n-\alpha_0|  |M''_n(\alpha_0) + R_n(\widehat
\alpha_n , \alpha_0)|$. Then,
\[
\frac{\sqrt{n}}{\sigma_1(\theta)}   |\widehat{\alpha
}_n-\alpha_0
 |
> 2 \frac{\sqrt{\ln n}}{m(\theta)} \Longrightarrow
\frac{\sqrt{n}}{\sigma_1(\theta)}  \bigl (
 |M'_n(\alpha_0) | + r_n \bigr)
> 2  \frac{\sqrt{\ln n}}{m(\theta)}
 |M''_n(\alpha_0) + R_n(\widehat\alpha_n , \alpha_0) |,
\]
provided that $M'_n(\widehat\alpha_n)\neq M'_n(\alpha_0)$. Next,
introducing the event $ \{2|M''_n(\alpha_0)+R_n(\widehat\alpha_n,\break
\alpha_0)| > m(\theta)  \}$ and its complement (which includes
the event $\{M'_n(\widehat\alpha_n) = M'_n(\alpha_0)$\}),
we obtain
\begin{eqnarray*}
&&\P_{\theta}\biggl \{ \frac{ \sqrt{n} }{ \tau(\theta)}   |
\widehat{\alpha}_n-\alpha_0  | > 2 \sqrt{\ln n}
 \biggr\}\\
  && \quad \leq \P_{\theta} \biggl\{ \frac{ \sqrt{n} }{
\sigma_1(\theta) }   \{|M'_n(\alpha_0)| + r_n\} > \sqrt{\ln n}
 \biggr\}  +  \P_{\theta} \{2 |M''_n(\alpha_0) + R_n(\widehat\alpha_n ,
\alpha_0)|   \leq m(\theta) \}.
\end{eqnarray*}
It is easily checked from (A4) and $r_n =\mathrm{o}(n^{-1/2})$ that
\[
\sup_{\theta\in\Theta}
\P_{\theta} \biggl\{\frac{\sqrt{n}}{\sigma_1(\theta)} \{
|M'_n(\alpha_0)| + r_n\}
> \sqrt{\ln n}\biggr  \} = \mathrm{O}\biggl(\frac{1}{\sqrt{n}} \biggr) + 2
\Gamma \biggl(-\sqrt{\ln n} + \frac{\sqrt{n}r_n}{\sigma_1(\theta
)} \biggr) =
\mathrm{O}\biggl(\frac{1}{\sqrt{n}} \biggr).
\]
Finally, to obtain the bound (\ref{Lemme6}), it remains to justify the
use of the following bound:
%
\begin{equation} \label{eq6.1}
\sup_{\theta\in\Theta} \P_{\theta} \{2 |M''_n(\alpha_0)
+R_n(\widehat\alpha_n , \alpha_0) | \leq
m(\theta) \}=\mathrm{O}(n^{-1/2}) + \gamma_n.
\end{equation}

Using elementary inequalities and assumption (A5),
\begin{eqnarray*}
 &&\P_{\theta} \{2 |M''_n(\alpha_0) + R_n(\widehat\alpha
_n ,
\alpha_0)| \leq m(\theta) \}  \\
&& \quad  \leq
\P_{\theta} \{|M''_n(\alpha_0) -m(\theta)| \geq
m(\theta)/4 \} + \P_{\theta} \{| R_n(\widehat\alpha_n ,
\alpha_0)| \geq m(\theta)/4 \}\\
&& \quad \leq\P_{\theta} \{|M''_n(\alpha_0) - m(\theta)| \geq
m(\theta)/4 \} + \P_{\theta} \{[| \widehat\alpha_n -
\alpha_0)|+ \omega_n]W_n \geq m(\theta)/4 \}\\
& & \quad =: P_{1,n,\theta} + P_{2,n,\theta}.
\end{eqnarray*}
It follows from (A4) that
$a:=\inf_{\theta\in\Theta} (m(\theta)/4\sigma_2(\theta)
)$ is
well defined and positive, and
%
\begin{equation} \label{P1ntheta}
\sup_{\theta\in\Theta} P_{1,n,\theta} \leq \mathrm{O}(n^{-1/2} ) +
2  \Gamma\bigl(-a\sqrt n\bigr) =
\mathrm{O}(n^{-/1/2} ).
\end{equation}
Now, let $d (\theta):= m(\theta)/4c_W$ with $c_W$ defined in (A5) and
notice that $d=\inf_{\theta'\in\Theta} d(\theta')/2$ in (A6). Use
the event
$ \{|\widehat{\alpha}_n - \alpha_0  |\leq d(\theta)-\omega
_n \}$ and
its complement to write
\begin{eqnarray*}
P_{2,n,\theta} &\leq& \P_{\theta} \biggl\{\frac{m(\theta)}{4} \le
[|\widehat{\alpha}_n - \alpha_0 |+ \omega_n]  W_n \leq W_n
d(\theta) \biggr\} +
\P_{\theta} \{ |\widehat{\alpha}_n - \alpha_0| > d(\theta) -
\omega_n  \} \\
&\leq& \sup_{\theta\in\Theta} \P_{\theta} \{c_W \le W_n
 \} + \sup_{\theta\in\Theta}\P_{\theta} \{|\widehat
{\alpha}_n - \alpha_0| >
d  \} =
\mathrm{O}(n^{-1/2}) + \gamma_n,
\end{eqnarray*}
from (A5)--(A6) and provided that $ \omega_n \leq d$. Therefore,
inequality (\ref{eq6.1}) holds true.
\end{pf}

Now, it remains to investigate the case $|u|< 2\sqrt{\ln n}$.

\begin{pro} \label{Prop2}
There exists a positive constant $C$ such that, for any $|u| < 2\sqrt
{\ln n}$,
%
\begin{equation} \label{Ineg_BE_inf}
\sup_{\theta\in\Theta} \biggl|\P_{\theta} \biggl\{\frac{\sqrt
n}{\tau(\theta)}  (\widehat{\alpha}_n-\alpha_0) \leq u \biggr\} -
\Gamma(u) \biggr| \leq C \biggl( \frac{1}{\sqrt n} + \sqrt{n}  r_n +
\omega_n + \gamma_n \biggr).
\end{equation}
\end{pro}

\begin{pf}
We just have to prove that (\ref{Ineg_BE_inf}) holds true for all
$n\ge n_0$, for some $n_0\in\N$. Let us introduce some sets and
derive their probability bounds:
\begin{itemize}
\item $E_{n,\theta}:= \{ \sqrt{n}  |\widehat{\alpha
}_n-\alpha_0| / \tau(\theta) \le2 \sqrt{\ln n}  \}$.
From (\ref{Lemme6}), $\sup_{\theta\in\Theta}\P_{\theta
}(E_{n,\theta}^c) = \mathrm{O}(n^{-1/2} + \gamma_n)$.
\item $A_{n}:=\{0\leq W_n \le c_W \}$ where the r.v. $W_n$ and the
constant $c_W$ are defined
in (A5). Then $\sup_{\theta\in\Theta}\P_{\theta}(A_{n}^c) = \mathrm{O}(n^{-1/2})$.
\item $D_{n,\theta}:= \{ 2 M''_n(\alpha_0) > m(\theta)
\}$. We have $\P_{\theta}\{ D_{n,\theta}^c \} \leq
\P_{\theta} \{ |M''_n(\alpha_0)- m(\theta)  |\geq
m(\theta)/2 \} \le\P_{\theta} \{ |M''_n(\alpha_0)-
m(\theta)  |\geq m(\theta)/4 \} $.
We know from (\ref{P1ntheta}) that
$\sup_{\theta\in\Theta}\P_{\theta}(D_{n,\theta}^c) = \mathrm{O}(n^{-1/2})$.
\end{itemize}
Then, we obtain from the previous estimates that the following set
\[ \label{def_Bntheta}
B_{n,\theta} := E_{n,\theta} \cap A_{n} \cap D_{n,\theta}
\]
is such that
%
\begin{equation}\label{Bntheta_o}
\sup_{\theta\in\Theta}\P_{\theta}(B_{n,\theta}^c)
\le \mathrm{O}(n^{-1/2}+\gamma_n).
\end{equation}

Now, if $D_{n,\theta,u} := \{\sqrt n
(\widehat{\alpha}_n - \alpha_0)/ \tau(\theta) \leq u\}$, then
we can write from (\ref{Bntheta_o})
%
\begin{equation}\label{etoile2}
 | \P_{\theta}(D_{n,\theta,u}) -\Gamma(u)  |
 \le | \P_{\theta}(D_{n,\theta,u}\cap B_{n,\theta} )
-\Gamma(u) | + \mathrm{O}(n^{-1/2}+\gamma_n).
\end{equation}
From (A2) and (A4), $ 0 <\overline{\sigma}:=\sup_{\theta\in\Theta
} \tau(\theta) < \infty$. Define the piecewise quadratic functions
\begin{equation}\label{g_plus_et_moins}
g^-(v) := c^- + b^- v + a^- v^2,  \qquad   g^+(v) := c^+ + b^+ v + a^+
v^2,
\end{equation}
 where   $c^{\pm}:= n [ M'_n(\alpha_0) \pm r_n],   b^{\pm
}:=\tau(\theta)\sqrt{n}  [ M''_n(\alpha_0) \pm \operatorname{sign}(v) c_W
\omega_n ],
  a^{\pm}:=   \pm \overline{\sigma}^2   c_W,$
and $\operatorname{sign} (v)$ denotes the sign of $v$ when $v\neq0$ and $\operatorname{sign}
(0)=0$. Notice that $g^-$ and $g^+$ are continuous on the whole real
line. To bound the term $ | \P_{\theta}(D_{n,\theta,u}\cap
B_{n,\theta} ) -\Gamma(u) |$ in (\ref{etoile2}), let us
introduce the events
%
\begin{equation} \label{def_Enthetat}
E_{n,\theta,u}^{\pm} := \{ g^{\pm}(u) \ge0\}.
\end{equation}
It follows from Lemma~\ref{prop2.1} in Appendix~\ref{B} that, for $n$
large enough
and $|u|< 2\sqrt{\ln n}$,
\[
\P_{\theta}(E_{n,\theta,u}^- \cap B_{n,\theta}) \le\P_{\theta
}(D_{n,\theta,u} \cap B_{n,\theta}) \le\P_{\theta}(E_{n,\theta
,u}^+ \cap B_{n,\theta})
\]
so that
\begin{eqnarray} \label{basic_estimate}
&& |\P_{\theta}( D_{n,\theta,u}\cap B_{n,\theta}) - \Gamma
(u) |\nonumber\\
&& \quad \leq\max \{  |\P_{\theta} (E^{-}_{n,\theta,u}\cap
B_{n,\theta} ) - \Gamma(u) | ,  |\P_{\theta} (
E^{+}_{n,\theta,u}\cap B_{n,\theta} ) - \Gamma(u) |
 \}\\
&& \quad \leq\max \{  |\P_{\theta} (E^{-}_{n,\theta,u} ) -
\Gamma(u) | ,  |\P_{\theta} ( E^{+}_{n,\theta,u} )
- \Gamma(u) |  \} + \P_{\theta}(B^c_{n,\theta}).
\nonumber
\end{eqnarray}
Then the proof of Proposition~\ref{Prop2} is easily completed using
(\ref{Bntheta_o}) and the following estimate:
There exists a constant $C$ such that for $n$ large enough and $|u|<
2\sqrt{\ln n}$
%
\begin{equation} \label{EstimateFinal}
\sup_{\theta\in\Theta} |\P_{\theta} (E_{n,\theta,u}^{\pm
} ) - \Gamma(u)  | \le C \biggl( \frac{1}{\sqrt n} +
\sqrt{n} r_n + \omega_n  \biggr).
\end{equation}

Indeed, $E^{\pm}_{n,\theta,u}=\{g^{\pm}(u)\ge0\}$ with
$g^{\pm}$ defined in (\ref{g_plus_et_moins}). We can write
\begin{eqnarray*}
E^{\pm}_{n,\theta,u} & = &  \bigl\{ n[ M'_n(\alpha_0) \pm r_n] + u
\tau(\theta) \sqrt{n} [ M''_n(\alpha_0) \pm \operatorname{sign}(u)c_W  \omega_n]
\pm u^2 \overline{\sigma}^2  c_W \ge0 \bigr\} \\
& = &  \biggl\{ \frac{\sqrt{n}}{\sigma_{n,u}(\theta)}
\biggl(M'_n(\alpha_0) + \frac{u  \sigma_1(\theta)}{\sqrt{n} m(\theta)}
 \bigl(M''_n(\alpha_0)- m(\theta) \bigr) \biggr) \geq - \frac
{a_n(u,\theta)+ b_n(u,\theta)}{\sigma_{n,u}(\theta)} \biggr\},
\end{eqnarray*}
where the positive real number $\sigma_{n,u}(\theta)$ is that of
condition (A4$'$) and
\[
a_n(u,\theta) = u   \biggl[\sigma_1(\theta) (1\pm\frac{\operatorname{sign}(u)c_W
\omega_n}{m(\theta)}  ) \pm\frac{u \overline{\sigma}^2
c_W}{\sqrt{n}}  \biggr],  \qquad   b_n(u,\theta) = \pm\sqrt{n} r_n .
\]
From the second statement of (A4$'$), it follows that there exists a
constant $B'$ such that we have, for $n$ large enough and $|u|<
2\sqrt{\ln n}$,
\[
\sup_{\theta\in\Theta} \biggl|\P_{\theta}(E_{n,\theta,u}^{\pm}) -
\Gamma \biggl(\frac{a_n(u,\theta)+ b_n(u,\theta)}{\sigma
_{n,u}(\theta)} \biggr)  \biggr| \le\frac{B'}{\sqrt n}.
\]
Now, from $\underline{\sigma}_1 :=
\inf_{\theta\in\Theta}\sigma_1(\theta) > 0$
and from the first property of $\sigma_{n,u}(\theta)$ in (A4$'$), it
follows that, for $n$ large enough and $|u|< 2\sqrt{\ln n}$,
and for all $\theta\in\Theta$, we have $\sigma_{n,u}(\theta) \geq
\underline{\sigma}_1/2$ and
\begin{eqnarray*}
 \biggl|\frac{a_n(u,\theta)}{\sigma_{n,u}(\theta)}-u \biggr|
&\leq& \frac{|u|}{\sigma_{n,u}(\theta)}\biggl ( |\sigma
_{n,u}(\theta) - \sigma_1(\theta) |
+\frac{c_W\omega_n}{m(\theta)}+
\frac{|u| \overline{\sigma}^2 c_W}{\sqrt{n}} \biggr) \\
&\leq& \frac{2|u| }{\underline{\sigma}_1} \biggl [  (A' +
\overline{\sigma}^2 c_W  )  \frac{|u|}{\sqrt
n} + \frac{c_W}{\underline m}\omega_n \biggr] \le
C'  \biggl(\frac{u^2}{\sqrt{n}}+|u|\omega_n \biggr),
\end{eqnarray*}
where $C'$ is independent of $n$, $u$, $\theta$.
We obtain from estimates on the characteristic function of the standard
Gaussian distribution reported in \cite{Pfa71}, page 89, that, for $n$
large enough, $|u|< 2\sqrt{\ln n}$, and
$\theta\in\Theta$,
\[
 \biggl|\Gamma \biggl(\frac{a_n(u,\theta)}{\sigma_{n,u}(\theta
)} \biggr) -\Gamma(u) \biggr| \leq C_1 \biggl(\frac{1}{\sqrt n}+\omega
_n \biggr)
\]
for some $C_1>0$.
We deduce from similar arguments that, for some constant $C_2$,
\[
 \biggl|\Gamma \biggl(\frac{a_n(u,\theta)}{\sigma_{n,u}(\theta
)} \biggr)
-\Gamma\biggl (\frac{a_n(u,\theta)+
b_n(u,\theta)}{\sigma_{n,u}(\theta)} \biggr)  \biggr| \le
C_2\sqrt{n}  r_n.
\]
Since $C_1$, $C_2$ only depend on $A'$,
$\underline\sigma_1$, $\underline m$, $\overline\sigma$ and $c_W$,
the proof of (\ref{EstimateFinal}) is complete.
\end{pf}

\section{A Berry--Esseen bound for an additive functional of
geometrically ergodic Markov chains} \label{IIbis}

The main focus of the paper is to apply the general Berry--Esseen
result of Theorem \ref{th-pfanzagl} to the case of $M$-estimators as
defined in the \hyperref[intro]{Introduction} when the observations come from a
geometrically ergodic Markov chain.
To check conditions (A4) and (A4$'$) in Theorem~\ref{th-pfanzagl}, we
need the next probabilistic results based on a recent version of the
Berry--Esseen theorem derived by \cite{HerPen08} in the geometrically
ergodic Markov chain setting.

\subsection{The statistical model} \label{model_stat}

\newcommand{\mub}{\overline{\mu}}

Let $(E,\cE)$ be a measurable space with a countably generated $\sigma
$-field $\cE$ and $\Theta$ be some general parameter space. Let $\{
X_n\}_{n\ge0}$ be a Markov chain with state space $E$, transition
kernels $\{Q_{\theta}(x,\cdot),   x\in E\}$, $\theta\in\Theta$
and an initial distribution $\mu$ that may or may not depend on
$\theta$.

\renewcommand{\theassumption}{($ \mathcal{M}$)}
\begin{assumption}\label{asumM}   Let $V \dvtx  E
{\rightarrow}[1, + \infty)$ be an unbounded function (independent of
$\theta$). For each $\theta\in\Theta$, there exists a $Q_\theta
$-invariant probability distribution, denoted by $\pi_\theta$, such that
\begin{enumerate}[(VG2)]
\item[(VG1)] $ b_1 := \sup_{\theta\in\Theta}  \pi
_\theta(V) < +\infty$.
\item[(VG2)] For all $\gamma\in(0,1]$, there exist real numbers
$\kappa_\gamma<1$ and $C_\gamma\geq0$ such that we have, for any
$\theta\in\Theta$, $n\geq1$ and $x\in E$,
\[
\sup \{ |Q_\theta^nf(x)-\pi_\theta(f) |, f  \dvtx
E{\rightarrow}\C\mbox{ measurable, }|f|\leq V^\gamma  \}\leq
C_\gamma  \kappa_\gamma^n  V(x)^\gamma.
\]
\end{enumerate}
\end{assumption}

Throughout Section \ref{IIbis}, we assume that $\mub(V):=\sup
_{\theta\in\Theta} \mu(V)<\infty$. Notice that (VG2) with $\gamma
= 1$ implies the following property: For any measurable real-valued
function $f$ defined on $E$ such that $|f|\leq D V,$ for some constant $D>0$,
%
\begin{equation}\label{sup_sim}
\forall n\geq1  \qquad  \sup_{\theta\in\Theta}  | \mathbb
{E}_{\theta, \mu} [f(X_n)] - \pi_\theta(f)  | \leq D C_1
\kappa_1^n   \mub(V).
\end{equation}
Moreover, conditions (VG1) and (VG2) imply that, for any $\gamma\in
(0,1]$ and $\theta\in\Theta$, $Q_\theta$ is $V^\gamma
$-geometrically ergodic, but it is worth noticing that the constants
$C_\gamma$ and $\kappa_\gamma$ do not depend on~$\theta$. In the
following remark, the properties (VG1) and (VG2) are related to the
so-called drift condition w.r.t. the function $V$ for each $Q_\theta$.

\begin{rem} \label{gamma=1}
Assume that for each $\theta\in\Theta$, $Q_{\theta}$ is aperiodic
and $\psi$-irreducible w.r.t. a~certain positive $\sigma$-finite
measure $\psi$ on $E$ (which may depend on $\theta$).
\begin{longlist}[2.]
\item[1.] For $\gamma=1$ and any fixed $\theta$, the properties
(VG1)--(VG2) follow from the drift condition: $Q_\theta V \leq
\varrho V +\varsigma  1_S$, with $\varrho<1,\varsigma>0$ and some
set{} $S$ ($S$ is the so-called small set) satisfying the minorization
condition $Q_\theta(x,\cdot) \geq c  \nu(\cdot)  1_S(x)$, where
$c>0$ and $\nu$ is a~probability measure concentrated on $S$ (see
\cite{MeyTwe93}, Theorem~16.0.1). In addition, the constants~$C_1$ and
$\kappa_1$ can be bounded by a quantity involving $\varrho$,
$\varsigma$, $c$, the measure $\nu$ and the set $S$ (see~\cite
{MeyTwe94}). To obtain the uniformity in $\theta$, it suffices to
check that all these elements do not depend on $\theta$.
\item[2.] \label{vg1-vg2} For any $\gamma\in(0,1]$, we have $\pi
_\theta(V^\gamma) \leq \pi_\theta(V)$ and thus condition (VG1)
implies that $\sup_{\theta\in\Theta}  \pi_\theta(V^\gamma) <
\infty$. Furthermore, under the drift condition, it follows from
Jensen's inequality that $Q_\theta V^\gamma\leq\varrho^\gamma V +
\varsigma^\gamma  1_S$. Using \cite{MeyTwe94}, one obtains (VG2).
\end{longlist}
\end{rem}

\subsection{A preliminary uniform Berry--Esseen statement}\label{II}

Let $\alpha_0 = \alpha_0 (\theta)\in\mathcal{A}$ be the parameter
of interest for the statistical applications we have in mind (see
condition (A1), page \pageref{conda1}), where $\theta$ is the
parameter of the Markov chain model and $\mathcal{A}$ is an open
interval of the real line.

Let $\xi(\alpha,x,y)$ be a real-valued measurable function
defined on $\mathcal{A}\times E^2$ such that the random variable $\xi
(\alpha,X_{k-1},X_k)$ is (integrable and) centered with respect to the
stationary distribution $\pi_\theta$, that is,
\[
\E_{\theta,\pi_\theta} [\xi(\alpha_0,X_{0},X_1)] = 0,
\]
and let
\[
S_n(\alpha) := \sum_{k=1}^{n} \xi(\alpha,X_{k-1},X_k).
\]
We investigate the following uniform Berry--Esseen property:
\[
\sup_{\theta\in\Theta}\sup_{u\in\R} \biggl|\P_{\theta,\mu
} \biggl\{\frac{S_n(\alpha_0)}{\sigma(\theta)\sqrt n} \leq u
\biggr\} - \Gamma(u) \biggr| = \mathrm{O}\biggl(\frac{1}{\sqrt n} \biggr),
\]
where $\sigma^2(\theta)$ will be defined below as the asymptotic
variance associated with the random variables $\xi(\alpha,X_{k-1},X_k)$.
When $\{X_n\}_{n\geq0}$ are i.i.d. and $\xi(\alpha,X_{k-1},X_k)
\equiv\xi(\alpha,X_k)$, this property follows from the Berry--Esseen
theorem \cite{Fel71}, provided that $\xi(\alpha,X_0)$ has finite
third-order moment, uniformly bounded in $\alpha$, and a variance
greater than some positive constant that does not depend on $\alpha$.

In our Markov framework, the following moment (or $V$-domination)
condition is natu\-ral for the functional $\xi$. In the sequel, this
condition will be required for $m_0 = 1,2$ or~3.
\begin{cnd*}   There exist real constants $m >
m_0\geq1$ and $C_\xi> 0$ such that
%
 \renewcommand{\theequation}{$D_{m_0}$}
 \begin{equation} \label{d-m}
\forall\alpha\in\mathcal{A}, \forall(x,y)\in E^2 \qquad    |\xi
(\alpha,x,y) |^m \leq C_\xi  \bigl (V(x) + V(y) \bigr).
\end{equation}
\renewcommand{\theequation}{\arabic{equation}}\setcounter{equation}{16}
\end{cnd*}

This domination condition implies that
%
\begin{eqnarray}\label{int_cond}
\E_{\theta,\pi_\theta}[  |\xi(\alpha,X_{0},X_1)|^m  ]  &=&
\int |\xi(\alpha,x,y)|^m Q_\theta(x,\mathrm{d}y)\,\mathrm{d}\pi_\theta(x)\nonumber
\\[-8pt]
\\[-8pt]
&\leq& C_\xi
  \bigl (\pi_\theta(V)  +  \pi_\theta(Q_\theta V) \bigr) <
\infty,
\nonumber
\end{eqnarray}
and since $m\geq1$, observe that
$\E_{\theta,\pi_\theta}[  |\xi(\theta,X_{0},X_1)|  ] < \infty$.

\begin{pro} \label{variance-asymptotic}
Suppose that Assumption \ref{asumM} holds true and that
$\xi$ is centered and satisfies condition $(D_1)$. Then, we have $\sup
_{\theta\in\Theta}\sup_{n\geq1}  |\mathbb{E}_{\theta,\mu
}[S_n(\alpha_0)] | < \infty$. In particular, for each $\theta\in
\Theta$, $\lim_n\mathbb{E}_{\theta,\mu} [S_n(\alpha
_0)/n ] = 0$.
If, in addition, $\xi$ satisfies condition $(D_2)$, then for each
$\theta\in\Theta$, the non-negative real number
\[
\sigma^2(\theta) := 
\lim_n \frac{\mathbb{E}_{\theta,\mu}[S_n(\alpha_0)^2]}{n}
\]
is well defined and does not depend on $\mu$. Furthermore, the
function $\sigma^2(\cdot)$ is bounded on~$\Theta$, and there exists
a positive constant $C$, only depending on $C_\xi$ and $\mub(V)$,
such that
\[
\forall\theta\in\Theta \qquad  \forall n\geq1,
 \biggl|\sigma^2(\theta) - \frac{\mathbb{E}_{\theta,\mu
}[S_n(\alpha_0)^2]}{n} \biggr| \leq\frac{C}{n}.
\]
\end{pro}

Now, we are ready to state our uniform Berry--Esseen statement for
$S_n(\alpha_0)$.
\begin{theo} \label{Th-B-E-proba}
Let us assume that:
\begin{longlist}[3.]
\item[1.] Assumption  \ref{asumM}  holds true;
\item[2.] the functional $\xi$ is centered and satisfies condition $(D_3)$;
\item[3.]$\sigma_0^2 := \inf_{\theta\in\Theta}\sigma
^2(\theta) >0$.
\end{longlist}
Then, there exists a constant $B(\xi)$ such that
\[
\forall n\geq1 \qquad  \sup_{\theta\in\Theta}\sup_{u\in\R}
\biggl|\P_{\theta,\mu} \biggl\{ \frac{S_n(\alpha_0)}{\sigma(\theta
)\sqrt n} \leq u \biggr\} - \Gamma(u) \biggr| \leq
\frac{B(\xi)}{\sqrt n}.
\]
Furthermore, the constant $B(\xi)$ depends on the functional $\xi$,
but only through $\sigma_0$ and the constant $C_\xi$ of condition $(D_3)$.
\end{theo}

The fact that we look for a Berry--Esseen bound with a constant
$B(\xi)$ independent of~$\theta$ is natural given our main purpose,
that is, to prove a uniform Berry--Esseen theorem for $M$-estimators.

There are several methods for deriving Berry--Esseen bound for the
functionals of Markov chains (see \cite{Bol82,Jen89}). But to prove
Proposition~\ref{variance-asymptotic} and
Theorem~\ref{Th-B-E-proba}, we use the weak spectral method developed in
\cite{HerPen08}.{}
(A Berry--Esseen theorem is established in \cite{Her08} for sequences
of the form $\{\xi(X_k)\}_{k\ge0}$ under
the conditions $\mu(V) < \infty$ and $|\xi|^3 \leq C  V$; however, the
case of sequences of the form $\{\xi(X_{k-1},X_k)\}_{k\ge0}$ is
not a direct corollary of this work since the Markov chain $\{
(X_{k-1},X_k)\}_{k\ge0}$ may not be
geometrically ergodic.) This method allows us to control the constant
$B(\xi)$ as a function of $C_\xi$ for checking assumption (A4$'$) of
Theorem \ref{th-pfanzagl} (see the arguments following equation (\ref
{B-E-final}) in Section \ref{III.2}). This follows from the next key
technical result. Although the proof of the Berry--Esseen theorem only
requires Taylor expansions up to the order $m_0$ and Condition (\ref
{d-m}) with $m_0=3$, for the purpose of possible further applications,
Lemma~\ref{l-L-r} below is stated for any $m_0\in\N^*$.

\begin{lem}\label{l-L-r} If $\xi$ is centered and satisfies Condition
 (\ref{d-m}) with $m_0\in\N^*$, then there exists $\beta>0$
such that
%
\begin{equation} \label{formule-Qn-lbda-L-R}
\forall\theta\in\Theta, \forall n\geq1, \forall t\in[-\beta
,\beta] \qquad
\E_{\theta,\mu} \bigl[\mathrm{e}^{\mathrm{i}tS_n(\alpha_0)} \bigr] = \lambda_\theta
(t)^n \bigl (1+L_\theta(t)\bigr) + r_{\theta,n}(t),
\end{equation}
where $\lambda_\theta(\cdot)$, $L_\theta(\cdot)$ and $r_{\theta
,n}(\cdot)$ are some $m_0$ times continuously differentiable
functions from $[-\beta,\beta]$ into $\C$ satisfying
$\lambda_\theta(0) = 1, \lambda_\theta'(0)=0, L_\theta(0) = 0$
and $r_{\theta,n}(0) = 0$. Furthermore, there exists $\rho\in(0,1)$
such that we have for $\ell=0,\ldots,m_0$:
\[
G_\ell:= \sup\bigl \{\rho^{-n}  \bigl|r_{\theta,n}^{(\ell)}(t)\bigr|,
|t|\leq\beta,  \theta\in\Theta, n\geq1 \bigr\}< \infty.
\]
Finally, the constants $\beta$, $\rho$, $G_\ell$ and the following
ones (for $\ell=0,\ldots,m_0$),
\begin{eqnarray*}
E_\ell&:=& \sup\bigl\{\bigl|\lambda_\theta^{(\ell)}(t)\bigr|, |t|\leq\beta,
\theta\in\Theta\bigr\} < \infty, \\
     F_\ell&: =& \sup\bigl\{\bigl|L_\theta
^{(\ell)}(t)\bigr|, |t|\leq\beta, \theta\in\Theta\bigr\} < \infty,
\end{eqnarray*}
depend on $\xi$, but only through the constant $C_\xi$ of Condition
 (\ref{d-m}).
\end{lem}

Lemma~\ref{l-L-r} is proved in Section~\ref{proof-lemma1}. The
definition of $L_\theta(t)$ and $r_{\theta,n}(t)$ (see (\ref
{form-L-resolovent}) and~(\ref{form-r-resolovent}))
shows that the constants $F_\ell$ and $G_\ell$ also depend on $\mub
(V)$ (see Remark~\ref{constantes-bis}). Now Lemma~\ref{l-L-r} allows
us to derive Proposition~\ref{variance-asymptotic} and Theorem~\ref
{Th-B-E-proba}.

\begin{pf*}{Proof of Proposition~\ref{variance-asymptotic}}
Assume that $\xi$ is centered and satisfies (\ref{d-m}) with $m_0\in
\N^*$. Proceeding as in (\ref{int_cond}) and using (\ref{sup_sim}),
(VG1) and $\mub(V)<\infty$, we obtain that
%
\begin{equation}\label{ca_ne_saute_pas aux_yeux}
\sup_{\theta\in\Theta} \sup_{k\geq1 }\mathbb{E}_{\theta,\mu
}[  |\xi(\alpha_0,X_{k-1},X_k)|^m  ] < \infty  \qquad \mbox{for
some $m > m_0$.}
\end{equation}
Now assume $m_0=1$, and let $\phi(t) := \E_{\theta,\mu
}[\mathrm{e}^{\mathrm{i}tS_n(\alpha_0)}]$, $t\in\R$. Then $\phi'(0) = \mathrm{i}  \E_{\theta
,\mu}[S_n(\alpha_0)]$, but Lemma~\ref{l-L-r} also gives $\phi'(0) =
L_\theta'(0) + r_{\theta,n}'(0)$. Hence $\sup_{\theta\in\Theta
}\sup_{n\geq1 } |\E_{\theta,\mu}[S_n(\alpha_0)]| \leq F_1 + G_1$.
Next, assume $m_0=2$. From (\ref{ca_ne_saute_pas aux_yeux})
we have
$\mathbb{E}_{\theta,\mu}[  S_n(\alpha_0)^2  ]< \infty$, and thus
we can write $\phi''(0) = -\E_{\theta,\mu}[S_n(\alpha_0)^2]$, and
$\phi''(0) = n\lambda_\theta''(0) + L_\theta''(0) + r_{\theta
,n}''(0)$ by Lemma~\ref{l-L-r}. Thus we obtain
$|\lambda_\theta''(0) + \mathbb{E}_{\theta,\mu}[S_n(\alpha
_0)^2]/n| \leq (|L_\theta''(0)|+ |r_{\theta,n}''(0)|)/n \leq(F_2 +
G_2)/n$. Set $\sigma^2(\theta) := - \lambda''_\theta(0)$. Then
$\sup_{\theta\in\Theta} \sigma^2(\theta) \leq E_2$ (by Lemma~\ref
{l-L-r}), and the proof is complete with $C:=F_2 + G_2$.
\end{pf*}

\begin{pf*}{Proof of Theorem~\ref{Th-B-E-proba}} Recall that
$\xi$ is centered and satisfies condition ($D_3$). To prove the
result, we use Lemma~\ref{l-L-r} with $m_0=3$ and we adapt the
arguments of the i.i.d. case.
Recall that $\sigma^2(\theta) = - \lambda''_\theta(0)$. According to
the classical Berry--Esseen inequality (see \cite{Fel71}), we
must prove that for some suitable positive constant $c$, $\sup_{\theta
\in\Theta} A_n(\theta) =\mathrm{O}(n^{-1/2}),$ where
\[
A_n(\theta) := \int_{-c\sqrt n}^{c\sqrt
n} \biggl|\frac{\E[\mathrm{e}^{\mathrm{i}t\fracc{S_n(\alpha_0)}{\sigma(\theta)\sqrt n}}]
- \mathrm{e}^{-\fraca{t^2}{2}}}{t} \biggr|\,\mathrm{d}t .
\]
For the moment, we just assume that $0< c \leq\beta\sigma_0$, where
$\beta$ is the real number in Lemma~\ref{l-L-r}. Notice
that $|t|\leq c $ implies $|t/\sigma(\theta)| \leq\beta$ for all
$\theta\in\Theta$. Using Lemma~\ref{l-L-r}, we have
\begin{eqnarray*}
A_n (\theta) & \hspace*{2pt}\leq& \int_{-c\sqrt n}^{c\sqrt n}
 \biggl|\frac{\lambda_\theta (\fracc{t}{\sigma(\theta)\sqrt
n} )^n - \mathrm{e}^{-\fraca{t^2}{2}}}{t} \biggr|\,\mathrm{d}t\\
&&{}+ \int_{-c\sqrt n}^{c\sqrt n}  \biggl|\lambda_\theta \biggl(\frac
{t}{\sigma(\theta)\sqrt
n} \biggr) \biggr|^n
 \biggl|\frac{L_\theta (\fracc{t}{\sigma(\theta)\sqrt n}
)}{t} \biggr|\,\mathrm{d}t \\
& &{}
+ \int_{-c\sqrt n}^{c\sqrt n}  \biggl|\frac{r_{\theta,n} (\fracc
{t}{\sigma(\theta)\sqrt n} )}{t} \biggr|\,\mathrm{d}t\\
& :=& I_n(\theta) +
J_n(\theta) + K_n(\theta).
\end{eqnarray*}
By a Taylor expansion, for all $\theta\in\Theta$ and $|v|\leq c$,
\[
 \biggl|\lambda_\theta \biggl(\frac{v}{\sigma(\theta)} \biggr) - 1 +
\frac{v^2}{2} \biggr| \leq\frac{E_3}{6\sigma_0^3}  |v|^3,
\]
where $E_3$ is defined in Lemma~\ref{l-L-r}. Hereafter, set
$c := \min\{\beta\sigma_0,3\sigma_0^3/2E_3,\sqrt2\}$. From the last
inequality, deduce that for any $|v|\leq c$
\[
 \biggl|\lambda_\theta \biggl(\frac{v}{\sigma(\theta)} \biggr)
\biggr| \leq1 - \frac{v^2}{2} + \frac{v^2}{4} \leq \mathrm{e}^{-\fraca{v^2}{4}}.
\]
Therefore, for any $t\in\R$ such that $|t|\leq c \sqrt n$,
%
\begin{equation} \label{Z}
 \biggl|\lambda_\theta \biggl(\frac{t}{\sigma(\theta)\sqrt n}
\biggr) \biggr|^n \leq \mathrm{e}^{-\fraca{t^2}{4}}.
\end{equation}
Let us write
\[
\lambda_\theta \biggl(\frac{t}{\sigma(\theta)\sqrt n} \biggr)^n -
\mathrm{e}^{-\fraca{t^2}{2}} =
 \biggl(\lambda \biggl(\frac{t}{\sigma(\theta)\sqrt n} \biggr) -
\mathrm{e}^{-\fracc{t^2}{2n}} \biggr)   \sum_{k=0}^{n-1}
\lambda_\theta \biggl(\frac{t}{\sigma(\theta)\sqrt n} \biggr)^{n-k-1}
\mathrm{e}^{\fracc{-kt^2}{2n}}.
\]
Notice that $|\lambda_\theta(t/\sigma(\theta)\sqrt n) -
\exp(-t^2/2n)| \leq (a+E_3/6\sigma_0^3 )|t/\sqrt n|^{3}$
if $a :=
\sup_{|v|\leq c}|\psi^{(3)}(v)|$ with $\psi(v) := 6\exp(-v^2/2)$.
Moreover,
\[
\sum_{k=0}^{n-1}  \biggl|\lambda_\theta \biggl(\frac{t}{\sigma
(\theta)\sqrt n} \biggr) \biggr|^{n-k-1} \mathrm{e}^{- \fraca{kt^2}{2n}}
\leq\sum_{k=0}^{n-1} \mathrm{e}^{- \fracc{t^2(n-k-1)}{4n}}  \mathrm{e}^{-
\fracc{kt^2}{4n}} \leq b  n  \mathrm{e}^{-\fraca{t^2}{4}},\vspace*{-2pt}
\]
where $b := \sup_{|v|\leq c}\exp(v^2/4)$. Hence
\[
 \biggl|\lambda_\theta \biggl(\frac{t}{\sigma(\theta)\sqrt n}
\biggr)^n - \mathrm{e}^{-\fraca{t^2}{2}} \biggr|
\leq \biggl(a+\frac{E_3}{6\sigma_0^3} \biggr)  b
n^{-\fraca{1}{2}}  |t|^{3}  \mathrm{e}^{-\fraca{t^2}{4}},\vspace*{-2pt}
\]
which yields $\sup_{\theta\in\Theta} I_n(\theta) \leq b  n^{-1/2}
 (a+ E_3/6\sigma_0^3 ) \int_\R t^{2} \exp(-t^2/4)\,\mathrm{d}t $.
Next, using (\ref{Z}) and $L_\theta(0) = 0$,
\[
\sup_{\theta\in\Theta} J_n(\theta) \leq\frac{  F_1}{\sigma
_0\sqrt
{n}}  \int_\R \mathrm{e}^{-\fraca{t^2}{4}}\,\mathrm{d}t.\vspace*{-2pt}
\]
Finally, using $r_{\theta,n}(0) = 0$, we have
$\sup_{\theta\in\Theta} |r_{\theta,n}(t/\sigma(\theta)\sqrt
n)| \leq(|t|/\sigma_0\sqrt n)  G_1  \rho^n$, so that
$\sup_{\theta\in\Theta} K_n(\theta) \leq(2 c G_1/\sigma_0)
\rho^n$. Gathering the results, we deduce that
\[
\sup_{\theta\in\Theta} A_n(\theta)\leq\frac{A}{\sqrt n} + \frac{2 c
G_1}{\sigma_0}  \rho^n,\vspace*{-2pt}
\]
where the constants $A, \rho,G_1$ and $c$
depend on $C_\xi$ of condition ($D_3$). The Berry--Esseen inequality
\cite{Fel71} then yields
\[
\sup_{u\in\R} \biggl|\P_{\theta,\mu}\biggl \{ \frac{S_n(\theta
)}{\sigma(\theta)\sqrt n} \leq u \biggr\} - \Gamma(u) \biggr| \leq
\frac{1}{\uppi} \biggl(\frac{A}{\sqrt n} + \frac{2 c G_1}{\sigma_0}
\rho^n + \frac{24\eta}{c\sqrt n} \biggr),\vspace*{-2pt}
\]
where $\eta=
\sup_{u\in\R}|\Gamma'(u)|$. The proof of Theorem~\ref
{Th-B-E-proba} is complete.\vspace*{-3pt}
\end{pf*}

\subsection{\texorpdfstring{Proof of Lemma~\protect\ref{l-L-r}}{Proof of Lemma 1}}\vspace*{-3pt} \label{proof-lemma1}

For $\theta\in\Theta$ fixed, Lemma~\ref{l-L-r} follows from \cite
{HerPen08}, Section~10. Here, we must prove that all the constants in
Lemma~\ref{l-L-r} are uniform in $\theta$ and depend on $\xi$ as
claimed. For this purpose, the weak spectral method is outlined below
(in the $V$-geometrical ergodicity context) and we give the main
statements by paying
special attention to the constants. For convenience, the technical
proofs are postponed in Appendix~\ref{A}.

$\bullet$ \textit{Geometrical ergodicity of $Q_\theta$.}
Let $0<\gamma\leq1$. We denote by $\cB_\gamma$ the weighted
supremum-normed space of measurable complex-valued functions $f$ on $E$
such that
\[
\|f\|_{\gamma} : = \sup_{x\in E} \frac{|f(x)|}{V(x)^\gamma} <
\infty.\vspace*{-2pt}
\]
$(\cB_\gamma,\|\cdot\|_\gamma)$ is a Banach space. The space of
bounded operators on $\cB_{\gamma}$ is denoted by $\cL(\cB_{\gamma
})$, and the associated operator norm is still denoted by $\|\cdot\|
_{\gamma}$. We have from (VG1)
%
\begin{equation} \label{pi-V-gamma}
\sup_{\theta\in\Theta}  \pi_\theta(V^\gamma) \leq b_1 = \sup
_{\theta\in\Theta}  \pi_\theta(V) < \infty,\vspace*{-2pt}\vadjust{\goodbreak}
\end{equation}
so that $\pi_\theta$ is a continuous linear form on $\cB_{\gamma}$.
Define the following rank-one projection on~$\cB_{\gamma}$:
\[
\forall f\in\cB_{\gamma} \qquad   \Pi_\theta f := \pi_\theta(f)1_E.
\]
Then condition (VG2) in Assumption \ref{asumM} can be rewritten as follows: $Q_\theta
\in\cL(\cB_{\gamma})$ and there exist $\kappa_\gamma<1$ and
$C_\gamma>0$ such that
%
\begin{equation} \label{VG2'}
\forall\theta\in\Theta,  \forall f\in\cB_{\gamma}, \forall
n\geq1 \qquad
\|Q_\theta^nf - \Pi_\theta f\|_{\gamma} \leq C_\gamma  \kappa
_\gamma^n  \|f\|_{\gamma}.
\end{equation}
From (\ref{pi-V-gamma}) and (\ref{VG2'}), $\|Q_\theta^n\|_\gamma=
\sup_{x\in E} (Q_\theta^n V^\gamma)(x)/V(x)^{\gamma}$ is uniformly
bounded in $n\in\N^*$ and $\theta\in\Theta$.

$\bullet$ \textit{The Fourier kernels associated with $Q_\theta
$ and $\xi$.}
Assume that, for all $\alpha\in\mathcal{A}$, $\xi(\alpha,\cdot
,\cdot)$ is measurable. The Fourier kernels associated with $Q_\theta
$ and $\xi$ are denoted by $\{Q_\theta(t)(x,\mathrm{d}y),t\in\R\}$ and
defined by
\[
\forall x\in E \qquad    Q_\theta(t)(x,\mathrm{d}y) := \mathrm{e}^{\mathrm{i} t \xi(\alpha
_0,x,y)}Q_\theta(x,\mathrm{d}y).
\]
Let us recall that $S_n(\alpha_0) := \sum_{k=1}^{n} \xi(\alpha
_0,X_{k-1},X_k)$. The following link between $Q_\theta(t)$ and the
characteristic function of $S_n(\alpha_0)$ is well-known in the
spectral method:
%
\begin{equation} \label{F}
\forall n\geq1,  \forall t\in\R \qquad \E_{\theta,\mu
}\bigl[\mathrm{e}^{\mathrm{i}tS_n(\alpha_0)}\bigr] = \mu (Q_\theta(t)^n1_E ).
\end{equation}
In fact, we have $\E_{\theta,\mu}[\mathrm{e}^{itS_n(\alpha_0)}  f(X_n)] =
\mu(Q_\theta(t)^nf)$ for any real-valued measurable bounded function
$f$ on $E$. This can be easily checked by induction using the Markov
property and the following equality:
\[
\forall n\geq2 \qquad  \E_{\theta,\mu}\bigl[\mathrm{e}^{\mathrm{i}tS_n(\alpha_0)}  f(X_n)\bigr]
= \E_{\theta,\mu} \bigl[\mathrm{e}^{\mathrm{i}tS_{n-1}(\alpha_0)}  (Q_\theta
(t)f)(X_{n-1}) \bigr].
\]

$\bullet$ \textit{Spectral study of $Q_\theta(t)$ on $\cB
_\gamma$ (for $t$ near 0).} It can be easily seen that, for all $t\in
\R$, we have $Q_\theta(t)\in\cL(\cB_{\gamma})$. For $\kappa\in
(0,1)$, we set
\[
\cD_\kappa:=\{z\in\C  \dvt  |z|\geq\kappa, |z-1| \geq(1-\kappa
)/2\}.
\]

\begin{lem} \label{lem-res-Qt}
Let $\gamma\in(0,1)$. For all $\kappa\in(\kappa_\gamma,1)$, there
exists $\beta_{\gamma,\kappa}>0$ such that, for $\theta\in\Theta
$, $|t|\leq\beta_{\gamma,\kappa}$ and $z\in\cD_\kappa$, we have
$ (z-Q_\theta(t) )^{-1}\in\cL(\cB_{\gamma})$ and
\[
\cR_{\gamma,\kappa} := \sup \bigl\{ \bigl\| \bigl(z-Q_\theta(t)
\bigr)^{-1} \bigr\|_{\gamma} \dvt \theta\in\Theta, |t|\leq\beta_{\gamma
,\kappa}, z\in\cD_\kappa \bigr\} <\infty.
\]
Moreover, the constants $\beta_{\gamma,\kappa}$ and $\cR_{\gamma
,\kappa}$ depend on $\xi$, but only via the constant $C_\xi$ of
Condition
 (\ref{d-m}).
\end{lem}

For $\theta$ fixed, Lemma~\ref{lem-res-Qt} is established in \cite
{HerPen08}, Proposition~10.1, thanks to the theorem of Keller and
Liverani \cite{KelLiv99,Liv04}. Here, we only have to prove that the
constants $\beta_{\gamma,\kappa}$ and~$\cR_{\gamma,\kappa}$ are
uniform in $\theta$ and depend on $\xi$ as stated above. According to
\cite{KelLiv99}, Remark, page~145, it is enough to check that the
constants are so involved in the hypotheses of the Keller--Liverani theorem.
This is due to Lemmas~\ref{resolvent-Q}--\ref{lem-K-L} in
Appendix~\ref{A}.

$\bullet$ \textit{Proof of formula} (\ref{formule-Qn-lbda-L-R}).
Now assume that $\xi$ satisfies Condition (\ref{d-m}) for some
\mbox{$m_0\in\N^*$}. Let $\gamma_0\in(0,1)$ be fixed such that $\gamma_0
+ m_0/m < 1$. For any $\kappa\in(\kappa_{\gamma_0},1)$, denote by
$\Gamma_{0,\kappa}$ the oriented circle centered at $z=0$, with
radius $\kappa$, and by $\Gamma_{1,\kappa}$ the oriented circle
centered at $z=1$, with radius $(1-\kappa)/2$. Note that both $\Gamma
_{0,\kappa}$ and $\Gamma_{1,\kappa}$ are contained in $\cD_\kappa
$. From (\ref{VG2'}) and Lemma~\ref{lem-res-Qt}, one can deduce that
we have, for all $n\geq1$, $\theta\in\Theta$, and $t\in[-\beta
_{\gamma_0,\kappa};\beta_{\gamma_0,\kappa}]$, the following
equality in $\cL(\cB_{\gamma_0})$:
%
\begin{equation} \label{Q_l_L_N}
Q_\theta(t)^n = \lambda_\theta(t)^n  \Pi_\theta(t)  +  N_\theta(t)^n,
\end{equation}
where $\lambda_\theta(t)$ is the dominating simple eigenvalue of
$Q_\theta(t)$ and $\Pi_\theta(t)$ and $N_\theta(t)^n$ are the
elements of $\cL(\cB_{\gamma_0})$ defined by the following line integrals:
\[
\Pi_\theta(t) := \frac{1}{2\mathrm{i}\uppi}\oint_{\Gamma_{1,\kappa}}
\bigl(z-Q_\theta(t) \bigr)^{-1}\,\mathrm{d}z    \quad \mbox{and} \quad     N_\theta(t)^n
:= \oint_{\Gamma_{0,\kappa}} z^n  \bigl (z-Q_\theta(t) \bigr)^{-1}\,\mathrm{d}z.
\]
Note that we have $\lambda_\theta(0) = 1$ and $\Pi_\theta(0) = \Pi
_\theta$ from (\ref{VG2'}). Also observe that, from Lemma~\ref
{lem-res-Qt} and the definition of $\Gamma_{0,\kappa}$, we have $\|
N_\theta(t)^n\|_\gamma= \mathrm{O}(\kappa^n)$. Since $1_E\in\cB_{\gamma
_0}$ and $\mu(V)<\infty$ ($\mu$~is a continuous linear form on $\cB
_{\gamma_0}$), the equalities (\ref{F}) and (\ref{Q_l_L_N}) give:
\[
\E_{\theta,\mu}\bigl[\mathrm{e}^{\mathrm{i}tS_n(\alpha_0)}\bigr] = \lambda_\theta(t)^n  \mu
 (\Pi_\theta(t)1_E ) + \mu ( N_\theta(t)^n1_E ).
\]
Therefore, formula (\ref{formule-Qn-lbda-L-R}) holds true with
\[
L_\theta(t) := \mu (\Pi_\theta(t)1_E ) - 1,  \qquad    r_{\theta
,n}(t) := \mu ( N_\theta(t)^n1_E ) \qquad   (n\in\N^*).
\]
We have $L_\theta(0) = \mu (\Pi_\theta1_E ) - 1 = 0$ and
$r_{\theta,n}(0) = \mu (N_\theta(0)^n1_E ) = \mu
(Q_\theta^n1_E - \Pi_\theta1_E ) = 0$. Finally, to make the link
with Lemma~\ref{regularity-eigen-element} below easier, let us observe that
\begin{eqnarray}\label{form-L-resolovent}
1+L_\theta(t) &=& \frac{1}{2\mathrm{i}\uppi}\oint_{\Gamma_{1,\kappa}} \mu
\bigl(\bigl(z-Q_\theta(t)\bigr)^{-1}1_E \bigr)\,\mathrm{d}z, \\\label{form-r-resolovent}
r_{\theta,n}(t) &=& \frac{1}{2\mathrm{i}\uppi} \oint_{\Gamma_{0,\kappa}} z^n
\mu \bigl(\bigl(z-Q_\theta(t)\bigr)^{-1}1_E \bigr)\,\mathrm{d}z.
\end{eqnarray}

$\bullet$ \textit{Regularity properties of $\lambda(\cdot)$,
$L_\theta(\cdot)$, $r_{\theta,n}(\cdot)$.} Let $\gamma_0'$ be such
that $\gamma_0 + m_0/m < \gamma_0' < 1$. We denote by $\cL(\cB
_{\gamma_0},\cB_{\gamma_0'})$ the space of the bounded linear
operators from $\cB_{\gamma_0}$ to $\cB_{\gamma_0'}$, and by $\|
\cdot\|_{\gamma_0,\gamma_0'}$ the associated operator norm.
\begin{lem} \label{regularity-eigen-element}
We have the following regularity properties:
\begin{longlist}[(a)]
\item[(a)] The map $Q_\theta(\cdot)$ is $m_0$-times continuously
differentiable from $\mathbb{R}$ to $\mathcal{L}(\mathcal{B}_{\gamma
_0},\mathcal{B}_{\gamma_0'})$, and we have
${\cal Q}_{\ell} : = \sup_{t\in\R,  \theta\in\Theta} \|Q_\theta
^{(\ell)}(t)\|_{\gamma_0,\gamma_0'} < \infty$ for $\ell=0,\ldots,m_0$.
\item[(b)] There exist some real numbers $\kappa\in(\kappa_{\gamma
_0},1)$ and $0<\beta<\beta_{\gamma_0,\kappa}$
such that, for all $\theta\in\Theta$ and $z\in\cD_\kappa$, the
function $R_{\theta,z} \dvtx  t\mapsto (z-Q_\theta(t) )^{-1}$ is
$m_0$-times continuously differentiable from $[-\beta,\beta]$ into
$\mathcal{L}(\mathcal{B}_{\gamma_0},\mathcal{B}_{\gamma_0'})$, and
we have for $\ell=0,\ldots,m_0$:
\[
\sup \bigl\{\bigl\|R_{\theta,z}^{(\ell)}(t)\bigr\|_{\gamma_0,\gamma_0'}
\dvt
|t|\leq\beta,  z\in\cD_\kappa, \theta\in\Theta \bigr\} <
\infty.
\]
\end{longlist}
The scalars $\beta$, $\kappa$ and all the bounds in \textup{(a)} and \textup{(b)}
depend on $\xi$ only via the constant~$C_\xi$ of Condition (\ref{d-m}).
\end{lem}

For $\theta$ fixed, Lemma~\ref{regularity-eigen-element} is
established in \cite{HerPen08}, Proposition~10.3. It can be
also derived from \cite{Gou08}, which relaxes the assumptions used in
\cite{HenHer04,GouLiv06} to obtain Taylor expansions of the resolvent
maps. (As observed in \cite{Gou08}, the passage to the
differentiability properties can be derived from \cite{Cam64}.)
However, a fine control of the constants is still required. Using
either~\cite{Gou08} or \cite{HerPen08}, Section~10, this control is
derived from Lemma~\ref{lem-res-Qt} and from Lemma~\ref{reg-Q(t)} in
Appendix~\ref{A}.

Since $1_E\in\cB_{\gamma_0}$ and $\mu$ is a continuous linear form
on $\cB_{\gamma_0'}$ (use $\mub(V)<\infty$), Lemma~\ref
{regularity-eigen-element}(b) gives that, for any $z\in\Gamma
_{0,\kappa}\cup\Gamma_{1,\kappa}$, the $\C$-valued function
$t\mapsto\mu((z-Q_\theta(t))^{-1}1_E)$ is $m_0$-times continuously
differentiable on $[-\beta,\beta]$ and that its $m_0$ first
derivatives are uniformly bounded in $\theta$ and $z\in\Gamma
_{0,\kappa}\cup\Gamma_{0,\kappa}$. The regularity properties (and
the related bounds) for~$L_\theta(\cdot)$ and $r_{\theta,n}(\cdot)$
then follow from (\ref{form-L-resolovent}) and (\ref
{form-r-resolovent}), while those concerning the function~$\lambda
_\theta(\cdot)$ follow from both Lemma~\ref
{regularity-eigen-element}(a) and Lemma~\ref
{regularity-eigen-element}(b), according to a formula given in~\cite
{HerPen08}, Section~7.2. Finally the property $\lambda_\theta'(0)=0$
can be proved as follows. By deriving~(\ref{formule-Qn-lbda-L-R})
(applied with $\mu=\pi_\theta$) at $t=0$ and by using the fact that
$\xi$ is centered, we have $0 = \mathrm{i}  \E_{\theta,\pi_\theta
}[S_n(\alpha_0)] = n  \lambda_\theta'(0) + L_\theta'(0) +
r_{\theta,n}'(0)$. Hence $\lambda_\theta'(0)=0$.
\begin{rem} \label{constantes-bis}
Notice that, according to (\ref{form-L-resolovent})--(\ref
{form-r-resolovent}), the constants $F_\ell$ and $G_\ell$ in
Lemma~\ref{l-L-r} also depend on the supremum in $\theta$ of the norm
of $\mu$ in $\mathcal{B}_{\gamma'_0}'$, namely $\sup_{\theta\in
\Theta}\mu(V^{\gamma'_0})$.
\end{rem}

\section{A Berry--Esseen theorem for $M$-estimators}
\label{III.2}

Consider a Markov chain satisfying Assumption \ref{asumM} of Section~\ref
{model_stat}. Let us introduce the statistic
%
\begin{equation} \label{mnf}
M_n(\alpha) := \frac{1}{n} \sum_{k=1}^n F(\alpha,X_{k-1},X_k),
\end{equation}
where $\alpha$ is the parameter of interest, $F$ is a real-valued
measurable function on $\mathcal{A}\times E^2$ and $\mathcal{A}$ is
an open interval of the real line.

Assume that
$F$ satisfies condition $(D_1)$ and let
\[
M_\theta(\alpha) :=
\lim_{n\rightarrow\infty}\E_{\theta,\mu} [M_n(\alpha)]= \E
_{\theta,\pi_\theta}[F(\alpha,X_0,X_1)] ,
\]
which is
well defined by Proposition~\ref{variance-asymptotic}. Assume also
that, for
each $\theta\in\Theta$, there exists a~unique $\alpha_0 =
\alpha_0(\theta)\in\mathcal{A}$, the so-called true value of the
parameter of interest, such that $M_\theta(\alpha) > M_\theta(\alpha
_0)$, $\forall\alpha\neq\alpha_0$.
To estimate $\alpha_0 = \alpha_0 (\theta),$ we consider an
$M$-estimator $\widehat\alpha_n$ as defined in Section \ref{IV.1},
that is,
$
M_n(\widehat\alpha_n )\leq\min_{\alpha\in\mathcal{A}} M_n(\alpha)
+c_n,
$ where $\{c_n\}_{n\geq1}$ is a sequence of non-negative real numbers
going to zero.

Let $F'$ and $F''$
be real-valued measurable functions defined on $\mathcal{A}\times E^2$
and let
%
\begin{equation}\label{mn_mn}
M'_n(\alpha) := \frac{1}{n} \sum_{k=1}^n F'(\alpha,X_{k-1},X_k), \qquad
  M''_n(\alpha) := \frac{1}{n} \sum_{k=1}^n
F''(\alpha,X_{k-1},X_k).
\end{equation}
The functionals $F'$ and $F''$ could be the first- and second-order
partial derivatives of~$F$ with respect to $\alpha$, but this is not
necessary to deduce our next result. Consider the following assumptions
on $F'$ and $F''$ (and, implicitly, on $c_n$; see (V3)).

\begin{assumptions*}
\begin{enumerate}[(V3)]
\item[(V0)] $ F'$ and $F''$ satisfy condition ($D_3$).
\item[(V1)] $ \forall\theta\in\Theta,  \E
_{\theta,\pi_\theta}[F'(\alpha_0,X_0,X_1)] = 0$ and $\alpha
_0=\alpha_0(\theta)$ is unique with this property.
\item[(V2)] $ m(\theta) := \E_{\theta,\pi_\theta}[F''(\alpha
_0,X_0,X_1)]$ satisfies
$ \inf_{\theta\in\Theta} m(\theta) >0 $.
\item[(V3)]
$M'_n(\widehat\alpha_n)$ satisfies condition \textup{(A3)}, that is,
$ \forall n\geq1$ and there exists $r_n >0$ independent of $\theta$
such that $r_n=\mathrm{o}(1/\sqrt{n})$ and
$\sup_{\theta\in\Theta} \P_{\theta,\mu}  \{|M'_n(\widehat
{\alpha}_n)|\geq r_n \} = \mathrm{O}(n^{-1/2}).$
\end{enumerate}
\end{assumptions*}

Notice that (V0) ensures $\sup_{\theta\in\Theta} m(\theta) <
\infty$ (see (\ref{int_cond})). Now, as a consequence of
Proposition~\ref{variance-asymptotic} applied to $F'$ and $F''$, the conditions
(V0)--(V2) enable us to define the asymptotic variances:
\begin{eqnarray*}
\sigma_1^2(\theta)& :=& \lim_n \frac{1}{n}
\mathbb{E}_{\theta,\mu} \Biggl[ \Biggl(\sum_{k=1}^n F'(\alpha
_0,X_{k-1},X_k) \Biggr)^2 \Biggr], \\
\sigma_2^2(\theta) & := & \lim_n \frac{1}{n}
\mathbb{E}_{\theta,\mu} \Biggl[ \Biggl(\sum_{k=1}^n F''(\alpha
_0,X_{k-1},X_k)- n  m(\theta) \Biggr)^2 \Biggr].
\end{eqnarray*}
Moreover, condition (V0) and Proposition
\ref{variance-asymptotic} ensure that $\sup_{\theta\in\Theta}
\sigma_j(\theta) < \infty$ for $j=1,2$. The following conditions are
also assumed to hold.
\textit{\begin{enumerate}[(V6)]
\item[(V4)] $  \inf_{\theta\in\Theta} \sigma_j(\theta) > 0$
 for  $j=1,2$.
\item[(V5)] There exist $\eta\in(0,1/2)$ and $C>0$ such that
\[
\forall(\alpha,\tilde\alpha)\in\mathcal{A}^2, \forall(x,y)\in
E^2 \qquad
 |F''(\alpha,x,y) - F''(\tilde\alpha,x,y) |\leq C  |\alpha
-\tilde\alpha|   \bigl(V(x)+ V(y) \bigr)^{\eta}.
\]
\item[(V6)] Set $d:=\inf_{\theta\in\Theta} m(\theta) /8\pi
_\theta(V^\eta)$ with $\eta$ defined in \textup{(V5)}. There exists
$\gamma_n=\mathrm{o}(1)$ such that
\[
\sup_{\theta\in\Theta}\P_{\theta,\mu} \{
|\widehat{\alpha}_n-\alpha_0| \geq d  \} \leq\gamma_n .
\]
\end{enumerate}
}

\begin{theo} \label{Th-B-E-estimator}
Assume that Assumption \ref{asumM} holds true, $F$ satisfies condition
$(D_1)$ and
conditions \textup{(V0)}--\textup{(V6)} are fulfilled. Let $\tau(\theta
) :=
\sigma_1(\theta)/m(\theta)$. Then there exists a positive constant
$C$ such that
\[
\forall n\geq1 \qquad  \sup_{\theta\in\Theta}\sup_{u\in\R}
\biggl|\P_{\theta,\mu} \biggl\{\frac{\sqrt n}{\tau(\theta)}  (\widehat
{\alpha}_n -\alpha_0) \leq u \biggr\}
- \Gamma(u) \biggr| \leq C  \biggl( \frac{1}{\sqrt n} +
\sqrt{n} r_n + \gamma_n  \biggr)  .
\]
\end{theo}

The statement in the above theorem corresponds to that of
the i.i.d. case in \cite{Pfa71} up to few
changes: First, the variances of the i.i.d. context{} (namely,
$\E_{\theta}[F'(\theta,X_0)^2]$ and
$\E_{\theta}[(F''(\theta,X_0) - m(\theta))^2]$ for an i.i.d. sequence
$\{X_n\}_{n\geq0}$ and a functional $F(\theta,x)) $ are replaced by
the above asymptotic variances $\sigma_1^2(\theta)$ and
$\sigma_2^2(\theta)$ (this is natural in a general
Markovian context); second, the uniform (in $\theta$) third-order
moment conditions{} (namely,
$\sup_{\theta\in\Theta}\E_{\theta}[  |F'(\theta,X_0)|^3 +
|F''(\theta,X_0)|^3  ] < \infty$) on both $F', F''$ are replaced by
the domination
condition ($D_3$) for $F', F''$; third, even when $F' = \partial
F/\partial\alpha$, here we allow for a~positive sequence $r_n$,
$n\geq1$, provided it decreases to zero sufficiently fast. The second
point is specific to the geometrically ergodic Markov chain case.
Indeed, in the same statistical model, Dehay and Yao \cite{DehYao07}
proved a CLT for maximum likelihood estimates under a~second-order
domination assumption on the two first derivatives of
the functional, which corresponds to inequality (\ref{d-m}) with
$m_0=2$. Here
the previous second-order assumption is replaced by the (almost)
optimal condition ($D_3$) for deriving the Berry--Esseen theorem for
$M$-estimators.

\begin{pf*}{Proof of Theorem~\ref{Th-B-E-estimator}} It suffices to check
the conditions (A1)--(A6) of Theorem~\ref{th-pfanzagl}. The limit
$M_\theta^\prime(\alpha):=\lim_n \E_{\theta,\mu} [M'_n(\alpha
)]$ is well defined by Proposition \ref{variance-asymptotic} and
condition (V0), the uniqueness of $\alpha_0$ is guaranteed by (V1) and
hence (A1) holds true. One more application of Proposition \ref
{variance-asymptotic} ensures that $\E_{\theta,\pi_\theta
}[F''(\alpha_0,X_0,X_1)]= \lim_n \E_{\theta,\mu} [M''_n(\alpha
_0)]$, hence (A2) is satisfied. Condition (V3) is nothing else but
(A3). The Berry--Esseen properties in (A4) are associated with the
functionals $F'(\alpha_0,x,y)$ and $F''(\alpha_0,x,y)$
respectively, so that they directly follow from Theorem~\ref{Th-B-E-proba}.

Now, let us check that (A5) holds true with $\omega_n\equiv0$.
Define $W := V^{\eta}$, where $\eta\in(0,1/2)$ is the scalar in (V5)
and notice that
$\E_{\theta,\pi_\theta}[W(X_0)^{1/\eta}] = \pi_\theta(V)$. Next,
since $V\ge1$
and $\eta\in(0,1/2)$, we have $1\le W \le W^2 \le V $ so that $1\le
\pi_{\theta}(W) \le\pi_{\theta}(W^2) \le\pi_{\theta}(V) \le b_1$
by property (VG1). Deduce that
$\sup_{\theta\in\Theta}\pi_{\theta}(W) < \infty$, and by Proposition
\ref{variance-asymptotic} applied to $\xi(\theta,x,y) = W(y)$
\[
\sup_{n\geq1}\sup_{\theta\in\Theta}\frac{1}{n}  \E_{\theta
,\mu} \Biggl[ \Biggl(\sum_{k=1}^n W(X_k)-n  \pi_{\theta}(W)
\Biggr)^2 \Biggr] < \infty.
\]
Now, condition (A5) is guaranteed by the properties \ref{asumM} and (V5)
with $\omega_n\equiv0$, $c_W := \sup_{\theta\in\Theta}
\pi_\theta(W)$ and $W_n := (1/n)  \sum_{k=1}^n(W(X_{k-1})
+W(X_k))$ provided that
%
\begin{equation} \label{eq6.2}
\sup_{\theta\in\Theta}\P_{\theta,\mu} \{ 8 \pi_{\theta}(W)
\le W_n
 \} = \mathrm{O}(n^{-1} ).
\end{equation}
To prove (\ref{eq6.2}), set $S_n := \sum_{k=1}^n W(X_k)$. Since $W_n
\leq 2  S_n/n +  (W(X_0) + W(X_n) )/n$ and
$\pi_{\theta}(W) \geq1$,
\begin{eqnarray*}
\P_{\theta,\mu} \{8 \pi_{\theta}(W) \le W_n  \} &\leq& \P
_{\theta,\mu} \{ S_n \geq 2  n  \pi_{\theta}(W) \} + \P
_{\theta,\mu} \{W(X_0)+W(X_n) \geq 4  n \pi_{\theta}(W) \}
\\
& \le& \P_{\theta,\mu} \{S_n - n\pi_{\theta}(W) \geq n
\} +
\P_{\theta,\mu} \{W(X_0)+W(X_n) \geq 4  n  \}.
\end{eqnarray*}
Equality (\ref{eq6.2}) is then obtained by Markov's inequality,
\begin{eqnarray*}
\P_{\theta,\mu} \{8 \pi_{\theta}(W)  \le  W_n  \}
&\leq&
\frac{1}{n^2}  \E_{\theta,\mu} \bigl[\bigl(S_n-n\pi_{\theta
}(W)\bigr)^2 \bigr] +
 \biggl(\frac{1}{4n} \biggr)^{1/\eta}  \E_{\theta,\mu}\bigl[\bigl(W(X_0) +
W(X_n)\bigr)^{1/\eta}\bigr]\\ & =& \mathrm{O}( n^{-1}  ),
\end{eqnarray*}
since
\[
\sup_{\theta\in\Theta}\sup_{n\geq1} \E_{\theta,\mu}\bigl[\bigl(W(X_0) +
W(X_n)\bigr)^{1/\eta}\bigr] \leq 2^{1/\eta-1}[ \mub(V) + C_1  \mub(V) +
b_1 ],
\]
using $(a +
b)^{1/\eta} \leq2^{1/\eta-1} (a^{1/\eta} +
b^{1/\eta} ) $ for any $a,b\geq0$ and (VG1)--(VG2). Notice also that
now condition (V6) is identical to condition (A6).

The difficult part is to check the Berry--Esseen-type property (A4$'$).
For this purpose, let $\Xi:= \{\xi_i(\cdot,\cdot,\cdot),
i\in I\}$ denote an arbitrary family of real-valued functionals
defined on $\mathcal{A}\times E^2$. Suppose that each $\xi_i$ is
centered, that is, $\E_{\theta,\pi_\theta}[\xi_i(\alpha
_0,X_0,X_1)] = 0$ for all
$i\in I$ and $\theta\in\Theta$, and that condition $(D_3)$ is
fulfilled uniformly in $i\in I$, that is,
%
\begin{equation} \label{D-Xi}
\exists m > 3, \exists C\geq0,  \forall i\in I, \forall
\alpha\in\mathcal{A}, \forall(x,y)\in E^2 \qquad
 |\xi_i(\alpha,x,y) |^m \leq C  \bigl (V(x) + V(y) \bigr).
\end{equation}
For each $i\in I$, set $S_n(\alpha_0,i) := \sum_{k=1}^{n}
\xi_i(\alpha_0,X_{k-1},X_k)$, and using Proposition
\ref{variance-asymptotic}, associate the corresponding asymptotic
variance denoted by $\sigma_i^2(\theta)$. Moreover, assume that
%
\begin{equation} \label{S-Xi}
0 < \inf\{\sigma_i(\theta), \theta\in\Theta, i\in I\} \leq
\sup\{\sigma_i(\theta), \theta\in\Theta, i\in I\} < \infty.
\end{equation}
Then, we deduce from Theorem~\ref{Th-B-E-proba} that, under
conditions \ref{asumM}, (\ref{D-Xi}), (\ref{S-Xi}) and $\mub(V) <
\infty$, there exists a constant $B$ such that
%
\begin{equation} \label{B-E-final}
\forall n\geq1 \qquad    \sup_{i\in I}\sup_{\theta\in\Theta}\sup
_{u\in
\R} \biggl|\P_{\theta,\mu} \biggl\{
\frac{S_n(\alpha_0,i)}{\sigma_i(\theta)\sqrt n} \leq u \biggr\} -
\Gamma(u) \biggr| \leq\frac{B}{\sqrt n}.
\end{equation}
This allows us to establish the two conditions in (A4$'$). Indeed, for
$(p,v)\in\N^*\times\mathbb{R}$ with~$v$
such that $|v|\leq2\sqrt{\ln p}$, let us introduce the functional
$\xi_{p,v}$ defined by
\[
\xi_{p,v}(\alpha_0,x,y) := F'(\alpha_0,x,y) + \frac{v}{\sqrt
p}\frac{\sigma_1(\theta)}{m(\theta)} \bigl(F''(\alpha_0,x,y) -
m(\theta) \bigr).
\]
Set $S_n(\alpha_0,p,v)  := \sum_{k=1}^{n}
\xi_{p,v}(\alpha_0,X_{k-1},X_k)$, and
\begin{eqnarray*}
\alpha_\theta(p,v)  &:= & \frac{v}{\sqrt p}\frac{\sigma_1(\theta
)}{m(\theta)}, \\
  S'_n(\theta)  &:=& \sum_{k=1}^n   F'(\alpha
_0,X_{k-1},X_k),\\
S''_n(\theta) &:=& \sum_{k=1}^n   F''(\alpha_0,X_{k-1},X_k) -n
m(\theta),
\end{eqnarray*}
so that $S_n(\alpha_0,p,v) = S'_n(\alpha_0) +
\alpha_\theta(p,v)  S''_n(\alpha_0)$. Notice that
$\E_{\theta,\pi_\theta}[\xi_{p,v}(\alpha_0,X_0,X_1)] = 0$ by~(V1)--(V2). We have
\begin{eqnarray*}
&&\mathbb{E}_{\theta,\pi_\theta} [ S_n(\alpha_0,p,v)^2 ] - \mathbb
{E}_{\theta,\pi_\theta}[ S'_n(\alpha_0)^2 ]\\
&& \quad  =
\alpha_\theta(p,v)^2 \mathbb{E}_{\theta,\pi_\theta}[ S''_n(\alpha
_0)^2  ] +
2  \alpha_\theta(p,v) \mathbb{E}_{\theta,\pi_\theta}[
S'_n(\alpha_0)
S''_n(\alpha_0) ].
\end{eqnarray*}
From (V2) and the fact that $\sigma_1(\cdot)$ is bounded, we have
$|\alpha_\theta(p,v)| \leq A
|v|/\sqrt p$ for some $A>0$ that does not
depend on $\theta$. Besides, as already mentioned in this
section, one can define the asymptotic variances
$\sigma_1^2(\theta)$ and $\sigma_2^2(\theta)$ associated with the
functionals $F'$ and~$F''$ by
\[
\sigma_1^2(\theta) := \lim_n \frac{1}{n}  \mathbb{E}_{\theta,\pi
_\theta}[  S'_n(\alpha_0)^2  ], \qquad
\sigma_2^2(\theta) := \lim_n \frac{1}{n}  \mathbb{E}_{\theta,\pi
_\theta}[
S''_n(\alpha_0)^2  ].
\]
Similarly, the asymptotic variance
$\sigma_{p,v}^2(\theta)$ associated with $\xi_{p,v}$ can be defined by:
\[
\sigma_{p,v}^2(\theta) := \lim_n \frac{1}{n}  \mathbb{E}_{\theta
,\pi_\theta}[  S_n(\alpha_0,p,v)^2].
\]
Then it follows from $ |  \mathbb{E}_{\theta,\pi_\theta}[
S'_n(\alpha_0)  S''_n(\alpha_0)   ] | \leq\mathbb{E}_{\theta
,\pi_\theta}[
S'_n(\alpha_0)^2  ]^{1/2}  \mathbb{E}_{\theta,\pi_\theta}[
S''_n(\alpha_0)^2   ]^{1/2}$ that
\[
 |\sigma_{p,v}^2(\theta) - \sigma_1^2(\theta) | \leq
A^2  \frac{v^2}{p} \sigma_2^2(\theta) +
2  A  \frac{|v|}{\sqrt p}  \sigma_1(\theta)  \sigma_2(\theta).
\]
Since $\sigma_j(\cdot)$ is bounded ($j=1,2$) and $|v|\leq2\sqrt{\ln
p} \leq2\sqrt p$, the previous inequality shows that there exists
$C'>0$, independent of $\theta$, such that
\[
 |\sigma_{p,v}^2(\theta) - \sigma_1^2(\theta) | \leq
C'  \frac{|v|}{\sqrt p} .
\]
Set $\overline{\sigma}_1 : = \sup_{\theta\in\Theta}
\sigma_1(\theta)$ and $\underline{\sigma}_1 : =
\inf_{\theta\in\Theta} \sigma_1(\theta)$ (we have
$\underline{\sigma}_1>0$ from (V4)). Using $|v|/\sqrt p\leq2
\sqrt{\ln p/p}$ and $\sqrt{\ln p/ p}=\mathrm{o}(1)$, the above inequality
implies that there exists
$P_0\in\N$ such that we have, for all $p\geq P_0$ and $v$ such that
$|v|\leq2\sqrt{\ln p},$
\[
\forall\theta\in\Theta \qquad   \tfrac{1}{2}  \underline{\sigma}_1
\leq\sigma_{p,v}(\theta) \leq\tfrac{3}{2}  \overline{\sigma}_1.
\]
In particular, under the same condition on $(p,v)$, this gives
$\sigma_{p,v}(\theta) + \sigma_1(\theta) \geq
3\underline{\sigma}_1/2$, hence $|\sigma_{p,v}(\theta) -
\sigma_1(\theta)| \leq2C' |v| / 3\underline{\sigma}_1\sqrt p$.
This proves the first assertion in (A4$'$).

Now, let us define
\[
I =  \bigl\{  (p,v)\in\N^*\times\mathbb{R} \dvt p\geq P_0, |v|\leq
2\sqrt{\ln p}   \bigr\}.
\]
It follows from (V0), (V2) and $\overline{\sigma}_1 < +\infty$ that
the family $\Xi: = \{\xi_{p,v},
(p,v)\in I\}$ satisfies~(\ref{D-Xi}). Besides, the above
bounds of $\sigma_{p,v}(\theta)$ give the property~(\ref{S-Xi}).
Then equation~(\ref{B-E-final}) shows that there
exists $B'>0$ such that we have for all $n\geq1$,
$(p,v)\in I$, $\theta\in\Theta$ and $u\in\R$:
\[
 \biggl|\P_{\theta,\mu} \biggl\{ \frac{S_n(\alpha_0,p,v)}{\sigma
_{p,v}(\theta)\sqrt n} \leq u \biggr\} - \Gamma(u) \biggr| \leq
\frac{B'}{\sqrt n}.
\]
Finally, let us fix any integer $n\geq P_0$ and
any real number $u$ such that $|u|\leq2\sqrt{\ln n}$. Then, the
previous Berry--Esseen bound with $p:=n$ and $v:=u$ provides the second
property of~(A4$'$). Indeed, we obtain\vadjust{\goodbreak}
from $S'_n(\alpha_0) = n  M'_n(\alpha_0)$ and $S''_n(\alpha_0) =
n  (M''_n(\alpha_0)-m(\theta))$ that
\begin{eqnarray*}
\frac{S_n(\alpha_0,n,u)}{\sigma_{n,u}(\theta)\sqrt n} &=& \frac
{1}{\sigma_{n,u}(\theta)\sqrt n}   \biggl(S'_n(\alpha_0) +
\frac{u}{\sqrt n}\frac{\sigma_1(\theta)}{m(\theta)}  S''_n(\alpha
_0) \biggr) \\
&=& \frac{\sqrt{n}}{\sigma_{n,u}(\theta)} \biggl(M'_n(\alpha_0) +
\frac{u \sigma_1(\theta)}{\sqrt{n} m(\theta)} \bigl(M''_n(\alpha
_0)-m(\theta) \bigr) \biggr).
\end{eqnarray*}
Now the proof of Theorem~\ref{Th-B-E-estimator} is complete.
\end{pf*}

\section{An example: AR(1) process with ARCH(1) errors} \label{iid}

Let us apply our theoretical results to an AR(1) process with ARCH(1)
errors that belongs to the class of ARMA--GARCH models (see \cite
{FraZak04} and the references therein). The observations are generated
by the process
%
%
\begin{equation}\label{ar_arch}
X_n = \rho_0 X_{n-1} + \sigma(X_{n-1}; a_0 , b_0 ) \varepsilon_n,
\qquad
  n=1,2,\ldots,
\end{equation}
where $X_0$ has some probability distribution $\mu$, $\sigma^2(x;a,b)
:= a + b x^2$ and $|\rho_0 |< 1$, $a_0,b_0 > 0$ are the true values of
the parameters. $\{ \varepsilon_n\}_{n\geq1}$ is a sequence of i.i.d.
random variables with zero mean and variance equal to 1, with finite
$p$th order moment for some $p$ to be specified below and (unknown)
density $f_\varepsilon$ that is continuous and positive on $\R$. $\{
\varepsilon_n\}_{n\geq1}$ is independent of $X_0$.
For simplicity, hereafter $\mu$ is assumed to be the Dirac
distribution~$\delta_0$. The ``true'' parameter $\theta$ in the
associated statistical model is the vector $(\rho_0, a_0, b_0)\in
\Theta\subset[-\overline\rho, \overline\rho]\times[m_a, M_a]
\times[m_b, M_b] \subset\R^3$, where $ \overline\rho\in(0,1)$,
$0< m_a < M_a < \infty$ and $0< m_b < M_b < 1$ are given such that
$ (\overline\rho +  \sqrt{M_b} )^p  \int_{\mathbb
{R}} (1 + |y|)^pf_{ \varepsilon}(y)\,\mathrm{d}y < 1$.
For illustration, we apply our results to estimate $\rho_0$ and $b_0$.

First, let us check that the Markov chain defined by (\ref{ar_arch})
satisfies Assumption \ref{asumM} of Section~\ref{model_stat} with $V(x) =
(1+|x|)^p$. To check (VG1)--(VG2) and the existence of the $Q_\theta
$-invariant probability measure $\pi_\theta$, by \cite{MeyTwe94},
Theorem~2.3, it suffices to prove that there exist constants $\varrho
\in(0,1)$, $c,\varsigma>0$, a Borel subset $S$ of the real line and a
probability measure~$\nu$ concentrated on $S$ such that the following
two conditions hold true (see Remark~\ref{gamma=1}):
For all $\theta\in\Theta$,
%
%
\begin{equation} \label{but-drift-mino}
\forall x\in\R \qquad    Q_\theta V(x)\leq\varrho  V(x) + \varsigma
  1_S(x)    \quad \mbox{and}  \quad   Q_\theta(x,\cdot) \geq c  \nu(\cdot
)  1_S(x).
\end{equation}
In our setting, the transition probability of $\{X_n\}_{n\geq0}$ is
given by
\[
Q_\theta(x,B) = \int1_B\bigl(\rho_0 x+\sigma(x,a_0,b_0)y\bigr) f_\varepsilon
(y)\,\mathrm{d}y
\]
for any Borel set $B\subset\R$. As a consequence, for all $\theta\in
\Theta$ and $x\in\R$,
\begin{eqnarray*}
\frac{Q_\theta V(x)}{V(x)}  &=&  \int_\R   \frac{V (\rho_0
x+\sigma(x,a_0,b_0)y )}{V(x)}  f_\varepsilon(y)\,\mathrm{d}y\\
&\leq & \int_\R   \biggl (\frac{1 + \overline\rho|x| + (\sqrt
{M_a} +  \sqrt{M_b}|x| ) |y|}{1+|x|} \biggr)^{ p }  f_\varepsilon
(y)\,\mathrm{d}y.
\end{eqnarray*}
By Fatou's lemma,
\[
\limsup_{|x|{\rightarrow}\infty} \biggl (\sup_{\theta\in
\Theta}\frac{Q_\theta V(x)}{V(x)} \biggr) \leq\bigl (\overline\rho+
\sqrt{M_b} \bigr)^p \int_{\mathbb{R}} (1+|y|)^pf_{\varepsilon}(y)\,\mathrm{d}y
=: \iota<1.
\]
Next, fix $\varrho\in(\iota,1)$. There exists $s>0$ such that for
each $|x|>s$, $ Q_\theta V(x)\leq\varrho  V(x)$ for all $\theta\in
\Theta$. Set $S:=[-s;s]$. For all $x\in S$ and $\theta\in\Theta$,
\[
Q_\theta V(x)   \leq  \int_\R \bigl(1 + \overline\rho s + \bigl(\sqrt
{M_a} + \sqrt{M_b}  s \bigr) |y| \bigr)^p  f_\varepsilon(y)\,\mathrm{d}y < \infty,
\]
so that the first condition in (\ref{but-drift-mino}) is guaranteed.
To check the second condition in (\ref{but-drift-mino}), define
\[
0<\delta(u) := \inf_{x\in S,  \theta\in\Theta}f_\varepsilon
\bigl(\sigma^{-1}(x,a_0,b_0) (u-\rho_0 x) \bigr),   \qquad    u\in\R.
\]
Then, for any $x\in S$, Borel set $B\subset\R$ and $\theta\in\Theta$,
\begin{eqnarray*}
Q_\theta(x,B)  &=&   \int_\R 1_B \bigl(\rho_0x+\sigma
(x,a_0,b_0)y \bigr)  f_\varepsilon(y)\,\mathrm{d}y\\
 &=&  \int_B   \frac
{f_\varepsilon   (\sigma^{-1}(x,a_0,b_0) (u - \rho_0
x) ) }{\sigma(x,a_0,b_0)}\,\mathrm{d}u \geq \int_B   \frac{\delta
(u)}{m_a}\,\mathrm{d}u.
\end{eqnarray*}
Define the measure $m(\mathrm{d}u) :=m_a^{-1} \delta(u)\,\mathrm{d}u$ and notice that
$m(S)>0$. We deduce from above that all $\theta\in\Theta$, $x\in S$
and Borel set $B\subset\R$,
\[
Q_\theta(x,B) \geq m(B) \geq m(B\cap S) = m(S)  \nu(B),
\]
where $\nu$ is the probability measure $\nu(B) := m(B\cap S)/m(S)$.
Hence the second condition in (\ref{but-drift-mino}) is fulfilled and
Assumption \ref{asumM} is satisfied for $\{X_n\}_{n\geq0}$ defined in
(\ref{ar_arch}).

Second, to estimate $\rho_0$, one can use the least-squares estimator,
\[
\widehat\rho_n := \frac{\sum_{k=1}^n X_k X_{k-1}}{\sum_{k=1}^n
X_{k-1}^2} = \arg\min_\rho\frac{1}{n} \sum_{k=1}^n F(\rho, X_{k-1},X_k),
\]
where $F(\rho, X_{k-1},X_k) := (X_{k}-\rho X_{k-1})^2$. We show that
the assumptions of Theorem~\ref{Th-B-E-estimator} are satisfied so
that we have a uniform Berry--Esseen bound for $\widehat\rho_n$. Fix
some $p>6$ and recall that $\int_\R|y|^p f_\varepsilon(y)\,\mathrm{d}y <\infty
$. Take $F^\prime(\rho, X_{k-1},X_k) := -2X_{k-1}(X_{k}-\rho
X_{k-1})$ and $F^{\prime\prime}(\rho, X_{k-1},X_k) := 2X_{k-1}^2$.
The conditions (V0) and (V1) are obviously fulfilled. Next, define
$m(\theta) := \E_{\theta,\pi_\theta}[F^{\prime\prime}(\rho_0,
X_{k-1},X_k)]$ and notice that $m(\theta)/2 = a_0 + (b_0+ \rho_0^2)
m(\theta)/2$. It follows that $m(\theta) = 2a_0/(1-\rho_0^2 - b_0)>
2m_a$ and thus (V2) holds. Condition (V3) is satisfied with $r_n \equiv0$.
From Proposition \ref{variance-asymptotic}, we can use the $Q_\theta
$-invariant probability measure $\pi_\theta$ to check condition (V4).
Notice that $\lim_n \E_{\theta,\pi_\theta}[X_{n}^2]=m(\theta)/2>
m_a$ and recall that $\{\varepsilon_n\}_{n\geq1}$ is i.i.d. We deduce that
\[
\sigma_1^2 (\theta) = \lim_n \frac{4}{n} \sum_{k=1}^n \E_{\theta
,\pi_\theta}  [ X_{k-1}^2\sigma^2(X_{k-1},a_0,b_0)\varepsilon
_k^2  ] \geq 4a_0 \lim_n \E_{\theta,\pi_\theta
}[X_{n}^2]\geq4m_a^2.
\]
To derive a lower bound for $\sigma_2^2 (\theta)$, let us decompose
\[
\E_{\theta,\pi_\theta}  \Biggl[\sum_{k=1}^n  \bigl( F^{\prime
\prime}(\rho_0, X_{k-1},X_k) - m(\theta)  \bigr)  \Biggr]^2 = \sum
_{k=1}^n v_{k,k} + 2\sum_{1\leq k<l\leq n} v_{k,l},
\]
where $v_{k,l} := \E_{\theta,\pi_\theta} [  ( F^{\prime
\prime}(\rho_0, X_{k-1},X_k) - m(\theta)  )  ( F^{\prime
\prime}(\rho_0, X_{l-1},X_l) - m(\theta)  ) ]$, $k\leq l$.
It is easily checked that $v_{k,l} = (\rho_0^2 + b_0)v_{k,l-1}$ for
$k<l$. In particular, this implies $v_{k,l}>0$, $k\leq l$. Next, by
elementary inequalities, we can obtain $ \inf_\theta\E_{\theta,\pi
_\theta}  [( F^{\prime\prime}(\rho_0, X_{0},X_1) - m(\theta))
^2 ] \ge K$ for some positive constant $K$ depending on the
variance of $\varepsilon^2_1$. Deduce that $\sigma_2^2 (\theta) \ge
K$, hence (V4) holds true. Condition (V5) is trivially satisfied. To
check the consistency of condition (V6), we take advantage of the
explicit form of $\widehat\rho_n$. Indeed, we have
\begin{eqnarray*}
\widehat\rho_n -\rho_0 &=& \frac{n^{-1}\sum_{k=1}^n  ( X_k
X_{k-1} - \rho_0 \E_{\theta,\pi_\theta} [X_1^2] ) - \rho_0
n^{-1}\sum_{k=1}^n  ( X_{k-1}^2 - \E_{\theta,\pi_\theta}
[X_1^2] )}{n^{-1}\sum_{k=1}^n  ( X_{k-1}^2 - \E_{\theta,\pi
_\theta} [X_1^2] ) + \E_{\theta,\pi_\theta} [X_1^2]}\\
&\hspace*{2.4pt}  =:&\frac{\Delta_{1n}-\rho_0 \Delta_{2n}}{\Delta_{2n}+\E
_{\theta,\pi_\theta} [X_1^2]}.
\end{eqnarray*}
By Chebyshev's inequality, for any $d>0$, $\P_{\theta,\delta_0 }\{
|\Delta_{1n}| >d\}\leq d^{-2}n^{-1}\E_{\theta,\delta_0}[n\Delta
_{1n}^2]$. Proposition~\ref{variance-asymptotic} guarantees that $\E
_{\theta,\delta_0}[n\Delta_{1n}^2]$ is uniformly bounded (with
respect to $\theta$). Similar arguments apply to $\Delta_{2n}$. Since
$\E_{\theta,\pi_\theta} [X_1^2]>m_a$ for all $\theta$, we deduce
that (V6) holds with $\gamma_n = \mathrm{O}(n^{-1})$. Finally, by Theorem \ref
{Th-B-E-estimator}, there exists $C>0$ such that
%
%
\begin{equation}\label{BE_ar_coeff}
\forall n\geq1 \qquad  \sup_{\theta\in\Theta}\sup_{u\in\R}
\biggl|\P_{\theta,\delta_0} \biggl\{\frac{\sqrt n}{\sigma_1(\theta
)m(\theta)^{-1}}  (\widehat{\rho}_n -\rho_0) \leq u \biggr\}
- \Gamma(u) \biggr| \leq\frac{C}{\sqrt n}.
\end{equation}

Third, let us now turn to the estimation of $b_0$. For this purpose,
assume that the $\varepsilon_n$'s have a moment of order $p$ for some
$p>12$. Recall that $a_0 = m(\theta) (1-\rho_0^2 - b_0)/2$ and notice
that $\tau_0^2 := m(\theta)/2$ is easily estimated by $\widehat\tau
_n^2 := n^{-1}\sum_{k=1}^n X_k^2$.
Next, define
\begin{eqnarray}
T_n(b; r, v ) := \frac{1}{n}\sum_{k=1}^n {\eta_k (b,r,v)}^2
\nonumber\\
\eqntext{\mbox{with } \eta_k(b,r,v) := (X_k - r X_{k-1})^2 - v (1 - r^{2} - b
) - b X_{k-1}^2 ,}
\\
\eqntext{\mbox{with } \displaystyle \frac{\partial T_n}{\partial b} (b; r, v ) = \frac
{2}{n} \sum_{k=1}^n (v-X_{k-1}^2) \eta_k (b,r,v), \frac
{\partial^2 T_n}{\partial b^2} (b; r, v ) = \frac{2}{n} \sum_{k=1}^n
(v-X_{k-1}^2)^2.}
\end{eqnarray}
If $\rho_0$ and $a_0$ were known, one could easily estimate $b_0$ by
least squares, more precisely by minimizing $T_n(b; \rho_0, \tau
_0^2)$ with respect to $b.$ With this idea in mind, our feasible
estimator of $b_0$ is defined as follows:
\[
\widehat b_n := \arg\min_{b\in[m_b, M_b]} M_n(b)   \qquad \mbox{with }
M_n(b) := T_n(b; \widehat\rho_n, \widehat\tau^2_n).
\]
Define $F^\prime(b,X_{k-1},X_k) := 2(\tau_0^2-X_{k-1}^2)  \eta_k
(b,\rho_0,\tau_0^2)$, $F^{\prime\prime} (b,X_{k-1},X_k) := 2(\tau
_0^2-X_{k-1}^2)^2$~and
$M^\prime_n(b) := \partial T_n/\partial b (b; \rho_0, \tau_0^2)$,
$M^{\prime\prime}_n(b) := \partial^2 T_n/\partial b^2 (b; \rho_0,
\tau_0^2 )$.\vadjust{\goodbreak} Let us point out that, in this case, $M^\prime_n(\cdot
)$ and $M^{\prime\prime}_n(\cdot)$ are only approximations of the
derivatives of $M_n(\cdot)$. Checking assumptions (V0)--(V2) is
obvious and therefore we skip the details. To check condition~(V3) for
$M^\prime_n(\widehat b_n)$, we use the decomposition $M^\prime
_n(\widehat b_n) = A_n + \Delta_{n} = A_n + \Delta_{1n} + \Delta
_{2n}+ \Delta_{3n}$ with
\begin{eqnarray*}
A_n &:=& \frac{2}{n}\sum_{k=1}^n(\tau_0^2  -  X^2_{k-1}) \eta_k
(\widehat b_n , \widehat\rho_n , \widehat\tau_n^2),\\[-2pt]
  \Delta_{n}&:=&\frac{2}{n}\sum_{k=1}^n(\tau_0^2  -  X^2_{k-1})
 \bigl( \eta_k (\widehat b_n , \rho_0 , \tau_0^2)  - \eta_k
(\widehat b_n , \widehat\rho_n , \widehat\tau_n^2)  \bigr),
\\[-2pt]
\Delta_{1n} &:=& \frac{4(\widehat\rho_n  - \rho_0)}{n}\sum
_{k=1}^n(\tau_0^2  -  X^2_{k-1}) (X_k  - \rho_0X_{k-1})X_{k-1}, \\[-2pt]
\Delta_{2n} &:=& -\frac{2(\widehat\rho_n -\rho_0)^2}{n}\sum
_{k=1}^n(\tau_0^2  -  X^2_{k-1}) X_{k-1}^2,\\[-2pt]
\Delta_{3n} &:=& 2 \{\widehat\tau_n^2(1-\widehat\rho_n^2
-\widehat b_n) - \tau_0^2(1-\rho_0^2 -\widehat b_n) \} (\tau_0^2
- \widehat\tau_n^2 +X_n^2/n).
\end{eqnarray*}
We check that each term satisfies condition (V3) with a suitable $r_n$.
First, we can write
\[
0=\frac{\partial M_n}{\partial b} (\widehat b_n) = A_n + B_n
 \qquad \mbox{with } B_n:=\frac{2(\widehat\tau_n^2- \tau_0^2)}{n}\sum
_{k=1}^n \eta_k (\widehat b_n , \widehat\rho_n , \widehat\tau_n^2).
\]
By elementary algebra $
B_n = 2(\widehat\tau_n^2- \tau_0^2)(\widehat b_n + \widehat\rho
_n^2)X_n^2/{n}.
$
Using the Berry--Esseen bound for~$\widehat\tau_n^2$ (see Theorem
\ref{Th-B-E-proba}) and Markov's inequality for $X_n^{2+a}$ for some
small $a>0$, we can prove that $\P_{\theta,\delta_0}\{|B_n|\geq
n^{-1}\} = \mathrm{O}(n^{-1/2})$ so that $\P_{\theta,\delta_0}\{|A_n|\geq
n^{-1}\} = \mathrm{O}(n^{-1/2})$.
By the bound in\vspace*{2pt} equation (\ref{BE_ar_coeff}), we have $\sup_\theta\P
_{\theta,\delta_0}\{|\widehat\rho_n -\rho_0|^j\geq n^{-j/2} \log
^{j/2} n\} = \mathrm{O}(n^{-1/2})$, $j=1,2$. Use this with $j=1$ and our Theorem
\ref{Th-B-E-proba} for the centered functional $\xi(X_k,X_{k-1}) =
(\tau_0^2 - X^2_{k-1}) (X_k - \rho_0X_{k-1})X_{k-1}$ to deduce that
$\P_{\theta,\delta_0} \{|\Delta_{1n}|\geq n^{-1}\log n\} =
\mathrm{O}(n^{-1/2})$. Next, the bound on $|\widehat\rho_n -\rho_0|^2$ and
Theorem~\ref{Th-B-E-proba} applied to the centered functional
$\xi(X_k,X_{k-1}) = (\tau_0^2 - X^2_{k-1}) X_{k-1}^2 - \tau_0^4 + \E
_{\theta,\pi_\theta} [X_{k-1}^4]$ allow us to deduce that $\P
_{\theta,\delta_0} \{|\Delta_{2n}|\geq n^{-1}\log n\} =
\mathrm{O}(n^{-1/2})$. Finally, use the Berry--Esseen bounds for $\widehat\rho
_n$ and $\widehat\tau_n^2$ and Markov's inequality for~$X_n^{2+a}$
with some $a>0$ to deduce that $\P_{\theta,\delta_0} \{|\Delta
_{3n}|\geq n^{-1}\log n\} = \mathrm{O}(n^{-1/2}).$ Combining these facts gives
that $M^\prime_n(\widehat b_n)$ satisfies condition (V3) with $r_n =
n^{-1}\log n$. Condition (V4) can be checked using similar arguments to
those used for $\widehat\rho_n$ and, therefore, the details are
omitted. Condition (V5) is trivially satisfied.
Finally, let us note that\looseness=-1
\[
\widehat b_n - b_0 = \frac{\sum_{k=1}^n(\widehat\tau_n^2 -
X_{k-1}^2) \eta_k(b_0,\widehat\rho_n, \widehat\tau_n^2)}{\sum
_{k=1}^n (\widehat\tau_n^2 - X_{k-1}^2)^2},
\]\looseness=0
and thus condition (V6) can be checked by arguments that we already
used in this example. We deduce from Theorem \ref{Th-B-E-estimator}
that, for some suitable $\tau(\theta)$,
\[
\forall n\geq1 \qquad  \sup_{\theta\in\Theta}\sup_{u\in\R}
\biggl|\P_{\theta,\delta_0} \biggl\{\frac{\sqrt n}{\tau(\theta)}
(\widehat{b}_n -b_0) \leq u \biggr\}
- \Gamma(u) \biggr| = \mathrm{O}\biggl(\frac{\log n}{\sqrt n}
\biggr).\vadjust{\goodbreak}
\]
The log factor in this Berry--Esseen bound is the price we pay for
estimating $b_0$ by a~simple two-step procedure, easy to implement,
where we first estimate $\widehat\rho_n$ and $\widehat\tau_n^2$ and
then we use the least-squares criterion $M_n(b)=T_n(b; \widehat\rho
_n,\widehat\tau_n^2)$. We feel that the log factor could be removed
by using a direct approach where the three parameters are estimated
simultaneously, but the investigation of this idea with Markov chain
data is left for future work.\looseness=-1

\section{Conclusion}
In this paper, we study the Berry--Esseen theorem for $M$-estimators
(or minimum contrast estimators) of some parameter $\alpha_0$ on the
real line. The estimators are defined from a~criterion based on a
functional $F(\alpha,X_{n-1},X_n)$ of the observation process $\{X_n\}
_{n\ge0}$. Our approach to
derive such bounds relies on Pfanzagl's method originally proposed for
i.i.d. observations \cite{Pfa71}. In a first step, Theorem~1 in \cite
{Pfa71} is extended to
obtain Berry--Esseen bounds for $M$-estimators based on any sequence of
observations satisfying suitable conditions. In a second step, the
specific case of $V$-geometrically ergodic Markov observations is
considered. We show that such Markov framework allows us to apply our
general result provided that $F$ and related functionals $F', F''$
satisfy suitable domination conditions. This result covers those
reported in \cite{Rao73,MilRau89}, which are proved under much
stronger moment conditions. We argue that the domination conditions
used in the present paper give an almost optimal treatment of
Berry--Esseen bounds for $V$-geometrically ergodic Markov chains. This
is possible due to the operator-type procedure developed in \cite{HerPen08}.

There are several possible extensions of our results. A straightforward
one is to follow the lines of the proof \cite{Pfa71}, Theorem~2, and
to consider an estimator of the standard deviation in the Berry--Esseen
bounds when this standard deviation depends on $\theta$ only through
$\alpha_0$. The details are omitted. Next, for more effective bounds,
we need to carefully evaluate the constants involved throughout the
paper. This is a direction of future work. Finally, there is no doubt
that the operator-type procedure in \cite{HerPen08} could be further
used in statistical applications with Markov models, in particular with
strongly ergodic Markov chains. This is under investigation.

\begin{appendix}

\section{\texorpdfstring{Complements for the proof of~Theorem~\protect\ref{th-pfanzagl}}
{Complements for the proof of Theorem 1}}
\label{B}

The reader is referred to Proposition~\ref{Prop2} and its proof for
the notation and the definitions used throughout this part. The
following lemma gives key properties of the random functions $g^{\pm}$.

\begin{alem} \label{basic_functions_g}
The following properties hold true.
\begin{longlist}[2.]
\item[1.] If $\nu_{n,\theta}:= \sqrt{n} (\widehat{\alpha}_n -
\alpha_0)/\tau(\theta)$, then $A_n \subset \{ g^-(\nu
_{n,\theta})\le0 \le g^+(\nu_{n,\theta}) \}$.
\item[2.] For $\omega\in D_{n,\theta}$, $g^{\pm}$ are increasing on
the interval $(-
2\sqrt{\ln n},2\sqrt{\ln n})$ provided that
%
\begin{equation}\label{cond_n_grand}
\sqrt{n}\geq\frac{2c_W }{\underline{m}}
\biggl [ \frac{4 \overline{\sigma}^2 \overline{m}   \sqrt{\ln
n}}{\underline{\sigma}_1} +\sqrt{n} \omega_n  \biggr].
\end{equation}
\end{longlist}
\end{alem}

\begin{pf}
We can write from assumptions (A5) and (A3)
\begin{eqnarray*}
 | n M'_n(\alpha_0) + (\widehat{\alpha}_n - \alpha_0) n
M''_n(\alpha_0)  |
& = &  | n M'_n(\widehat\alpha_n) - (\widehat\alpha_n - \alpha
_0) n R_n(\alpha_0, \widehat\alpha_n)  | \\
&\le& nr_n + n |\widehat\alpha_n - \alpha_0 |
|R_n(\alpha_0, \widehat\alpha_n)  | \\
& \le& nr_n + n |\widehat\alpha_n - \alpha_0 |  [
 |\widehat\alpha_n - \alpha_0 | + \omega_n  ] W_n.
\end{eqnarray*}
If $\omega\in A_n$, then
\[
 | n M'_n(\alpha_0) + (\widehat{\alpha}_n - \alpha_0) n
M''_n(\alpha_0)  |\le n  |\widehat\alpha_n -
\alpha_0 |^2 c_W + n \omega_n |\widehat\alpha_n -
\alpha_0 | c_W + nr_n.
\]
This last inequality is rewritten as
\[
n[ M'_n(\alpha_0)-r_n] + \tau(\theta) \sqrt{n}
 [ M''_n(\alpha_0) - \operatorname{sign}(\nu_{n,\theta}) c_W  \omega_n
]\nu_{n,\theta}-
\tau(\theta)^2 c_W \nu_{n,\theta}^2 \le0
\]
  and
  \[
n[ M'_n(\alpha_0)+r_n] + \tau(\theta)
\sqrt{n}  [ M''_n(\alpha_0) + \operatorname{sign}(\nu_{n,\theta}) c_W
\omega_n
 ]\nu_{n,\theta}+
\tau(\theta)^2 c_W \nu_{n,\theta}^2
\ge0 ,
\]
with $\nu_{n,\theta}:= \sqrt{n} (\widehat{\alpha}_n -
\alpha_0)/\tau(\theta)$. Since $0< \tau(\theta) \le
\overline{\sigma}$, we obtain that
\[
g^-(\nu_{n,\theta})\le0  \mbox{ and }   g^+(\nu_{n,\theta
}) \ge0.
\]

The second statement is proved as follows for $g^+$. Note
that $a^+>0$ and $g^+$ is
continuous. If we restrict $v<0$, the minimum of this quadratic
function $g^+(v)$ is achieved at
\[
v_{\min} = -\frac{b^+}{2a^+} = -\frac{\tau(\theta) \sqrt{n} [
M''_n(\alpha_0) - c_W \omega_n] }{2 \overline{\sigma}^2
c_W },
\]
or at the origin if $v_{\min} \geq0$. Now, if $\omega\in
D_{n,\theta}$ and $n$ satisfies condition (\ref{cond_n_grand}), it
is easy to check that
\[ \label{u_min}
v_{\min} < -2 \sqrt{\ln n}
\]
and $g^+$ is strictly increasing on $(0,\infty)$. Hence,
$g^+$ is increasing on $(- 2\sqrt{\ln n},2\sqrt{\ln n})$. Similar
arguments apply for $g^-$.
\end{pf}

\begin{alem} \label{prop2.1}
We have for $n$ large enough and $|u|<2\sqrt{\ln n}$
%
\begin{equation}\label{inclusions}
E_{n,\theta,u}^- \cap B_{n,\theta} \subset D_{n,\theta,u} \cap
B_{n,\theta} \subset E_{n,\theta,u}^+ \cap
B_{n,\theta}.
\end{equation}
\end{alem}

\begin{pf} It is understood below that $\omega\in B_{n,\theta}$.
Since $B_{n,\theta} \subset E_{n,\theta}\cap D_{n,\theta}$ and $|u|
< 2\sqrt{\ln n}$, the second statement in
Lemma~\ref{basic_functions_g} guarantees that for $n$ large enough
\[
\sqrt{n} (\widehat{\alpha}_n - \alpha_0)/\tau(\theta) \le u
 \quad \Longrightarrow \quad  g^+ \bigl(\sqrt{n} (\widehat{\alpha}_n - \alpha
_0)/\tau(\theta) \bigr) \le g^+(u).
\]
Since $B_{n,\theta} \subset A_{n}$, the first assertion in
Lemma~\ref{basic_functions_g} yields $g^+ (\sqrt{n}
(\widehat{\alpha}_n - \alpha_0)/\tau(\theta) ) \ge0$ so that
$g^+(u) \ge0$ when $\sqrt{n}
(\widehat{\alpha}_n - \alpha_0)/\tau(\theta) \le u$. This proves
the second inclusion in (\ref{inclusions}).

Next, assume that $g^-(u)\ge0$. Since $g^-$ is
increasing, we have
\[
\sqrt{n} (\widehat{\alpha}_n -
\alpha_0)/\tau(\theta) >u  \quad \Longrightarrow \quad  g^- \bigl(\sqrt{n}
(\widehat{\alpha}_n - \alpha_0)/\tau(\theta) \bigr) > g^-(u) \ge0.
\]
Since $B_{n,\theta}\subset A_{n}$, we know from
Lemma~\ref{basic_functions_g} that $g^-(\sqrt{n} (\widehat{\alpha}_n
- \alpha_0)/\tau(\theta)) \le0$ which is in contradiction with
the above inequality. Thus, $g^-(u)\ge0$ gives $\sqrt{n}
(\widehat{\alpha}_n - \alpha_0)/\tau(\theta) \leq u$.
\end{pf}

\section{\texorpdfstring{Complements for the proof of Lemma~\protect\ref{l-L-r}}
{Complements for the proof of Lemma 1}} \label{A}
A first step to control the constants in Lemma~\ref{lem-res-Qt} is to
study the resolvent map $(z-Q_\theta)^{-1}$ of the transition kernel
$Q_\theta$ acting on $\cB_{\gamma}$.
\begin{alem} \label{resolvent-Q} Let $\delta, r$ be such that $\kappa
_\gamma<r < 1$ and $0<\delta<1-r$. Then, for any $z\in\C$ such that
$|z|>r$ and $|z-1| > \delta$, the operator $z-Q_\theta$ is invertible
on $\cB_{\gamma}$, and we have:
\[
H_\gamma(\delta,r) := \sup \{\|(z-Q_\theta)^{-1}\|_{\gamma},
\theta\in\Theta, |z|>r, |z-1| > \delta \} < \infty.
\]
\end{alem}

\begin{pf} Let $g\in\cB_{\gamma}$, and let us write
$h_\theta= g - \pi_\theta(g)  1_E$. Since $\pi_\theta(h_\theta)
= 0$, it follows from (VG2) that
$\|Q_\theta^n h_\theta\|_{\gamma} \leq C_\gamma  \kappa_\gamma
^n  \|h_\theta\|_{\gamma}$. Now assume $|z|>r$. Then
\[
\sum_{k\geq0} |z|^{-(k+1)}  \|Q_\theta^k h_\theta\|_{\gamma} \leq
\frac{C_\gamma}{\kappa_\gamma}\sum_{k\geq0} \biggl(\frac{\kappa
_\gamma}{r} \biggr)^{k+1}  \|h_\theta\|_{\gamma} \leq \frac
{C_\gamma}{r-\kappa_\gamma}  \|h_\theta\|_{\gamma}.
\]
Thus, $\psi_\theta:= \sum_{k\geq0} z^{-(k+1)}  Q_\theta^k
h_\theta$ is absolutely convergent in $\cB_{\gamma}$, we have
$(z-Q_\theta)\psi_\theta= h_\theta$ and $\|\psi_\theta\|_{\gamma
} \leq C_\gamma  \|h_\theta\|_{\gamma}/(r-\kappa_\gamma)$.
Besides, if $z\neq1$, then we clearly have
\[
(z-Q_\theta) \biggl(\frac{\pi_\theta(g)}{z-1}  1_E \biggr) = \pi
_\theta(g)  1_E.
\]
Now assume $|z|>r$ and $|z-1| > \delta$. Then the function $f_\theta
:= (\pi_\theta(g)/(z-1))  1_E + \psi_\theta$ is such that
$(z-Q_\theta)f_\theta= g$. Thus $(z-Q_\theta)^{-1}g = f_\theta$.
From (\ref{pi-V-gamma}), we obtain $|\pi_\theta(g)| \leq\pi_\theta
(|g|) \leq \pi_\theta(V^\gamma)  \|g\|_{\gamma} \leq b_1  \|g\|
_{\gamma}$ and $\|h_\theta\|_{\gamma} = \|g - \pi_\theta(g)  1_E\|
_{\gamma} \leq (1+b_1 )\|g\|_{\gamma}$.
This gives: $\|f_\theta\|_{\gamma} \leq(b_1/\delta) \|g\|_{\gamma} +
C_\gamma(1+b_1) \|g\|_{\gamma}/(r-\kappa_\gamma)$, hence \mbox{$H_\gamma
(\delta,r) \leq[ b_1/\delta + C_\gamma(1+b_1)/(r-\kappa_\gamma) ]
< \infty$}.\\\mbox{}
\end{pf}

Second, the constants involved in the Doeblin--Fortet inequality and
the weak continuity condition of the Keller--Liverani theorem are
proved to be uniform in $\theta$ and to depend on $\xi$ only via the
constant $C_\xi$ of (\ref{d-m}). We appeal to \cite{KelLiv99},
Remark, page~145, and to the improvements given in \cite{Liv04}. In the
context of strongly ergodic Markov chains, the hypotheses resulting
from \cite{KelLiv99,Liv04} are stated in \cite{HerPen08}, Section~4,
and used here with the auxiliary norm $\|f\|_1 := \sup|f|/V$ on $\cB
_{\gamma}$. In the sequel, for $0<\gamma< \gamma' \leq1$, we denote
by $\cL(\cB_{\gamma},\cB_{\gamma'})$ the space of the bounded
linear operators from $\cB_{\gamma}$ to $\cB_{\gamma'}$, and by $\|
\cdot\|_{\gamma,\gamma'}$ the associated operator norm (with the
convention $\|\cdot\|_{\gamma} = \|\cdot\|_{\gamma,\gamma}$ when
$\gamma'=\gamma$).
\begin{alem}\label{lem-K-L}
Let $\gamma\in(0,1)$. We have:
\begin{longlist}[(b)]
\item[(a)]$\forall\theta\in\Theta, \forall t\in\R, \forall n\geq
1, \forall f\in\cB_{\gamma},
\|Q_\theta(t)^nf\|_{\gamma} \leq C_\gamma  \kappa_\gamma^n  \|f\|
_{\gamma} + b_1  \|f\|_1$;
\item[(b)]$\forall\theta\in\Theta, \forall t\in\R, \|Q_\theta(t) -
Q_\theta\|_{\gamma,1} \leq2^{2-\gamma}  {C_\xi}^{\fracb{1-\gamma
}{m}}  (E_\gamma+E_1)  |t|^{1-\gamma}  \|f\|_{\gamma}$,
\end{longlist}
where $E_\gamma:= \sup_{\theta\in\Theta}\|Q_\theta\|_\gamma, E_1
:= \sup_{\theta\in\Theta}\|Q_\theta\|_1$ and $C_\gamma$, $\kappa
_\gamma$, $b_1$ are defined in (\ref{pi-V-gamma}) and~(\ref{VG2'}).
\end{alem}

\begin{pf}
By using the inequality $\|Q_\theta(t)^nf\|_{\gamma} \leq\|
Q_\theta^n|f|  \|_{\gamma}$, assertion (a) easily follows from~(\ref
{VG2'}) and (\ref{pi-V-gamma}). To establish (b), let us recall that
we have from (\ref{d-m}) (use $V\geq1$)
\begin{eqnarray*}
|\xi(\theta,x,y)|^{1-\gamma} &\leq& C_\xi^{(1-\gamma)/m}
\bigl(V(x) + V(y) \bigr)^{1-\gamma} \\
&\leq& 2^{1-\gamma}  C_\xi^{(1-\gamma
)/m}  \bigl (V(x)^{1-\gamma} + V(y)^{1-\gamma} \bigr).
\end{eqnarray*}
Let $f\in\cB_{\gamma}$. From the definition of $Q_\theta(t)f $ and
the inequalities $|f| \leq V^\gamma\|f\|_{\gamma}$, $|\mathrm{e}^{\mathrm{i}a}-1|\leq
2|a|^{1-\gamma}$, we obtain that
\begin{eqnarray*}
 |(Q_\theta(t)f)(x) - (Q_\theta f)(x) | &\leq& \|f\|_{\gamma
}  \int_{E}  \bigl|\mathrm{e}^{\mathrm{i} t \xi(\alpha_0,x,y)} - 1 \bigr|
V(y)^\gamma
Q_\theta(x,\mathrm{d}y) \\
&\leq& 2^{2-\gamma}  {C_\xi}^{\fracb{1-\gamma}{m}}  |t|^{1-\gamma
}  \|f\|_{\gamma}  [V(x)^{1-\gamma}   (Q_\theta V^\gamma
 )(x) + (Q_\theta V)(x) ],
\end{eqnarray*}
from which we deduce (b).
\end{pf}

For the next lemma (used to prove Lemma~\ref
{regularity-eigen-element}), we introduce the following notation. For
any $\theta\in\Theta$, $k\in\N$, $t\in\mathbb{R}$, let us denote
by $Q_{\theta,k}(t)$ the operator associated with the kernel:
$Q_{\theta,k}(t)(x,\mathrm{d}y) = \mathrm{i}^k\xi(\alpha_0,x,y)^k \mathrm{e}^{\mathrm{i}t\xi(\alpha
_0,x,y)}  Q_\theta(x,\mathrm{d}y) $ ($x\in E$).
\begin{alem} \label{reg-Q(t)}
Let $0<\gamma< \gamma' \leq1$ and $k=0,\ldots,m_0$:
\begin{longlist}[(b)]
\item[(a)] If $\gamma+ k/m < \gamma'\leq1$, then the map $t\mapsto
Q_{\theta,k}(t)$ is continuous from $\mathbb{R}$ to $\mathcal
{L}(\mathcal{B}_{\gamma},\mathcal{B}_{\gamma'})$.
\item[(b)] If $k\leq m_0-1$ and $\gamma+(k+1)/m < \gamma'\leq1$, then the
map $t\mapsto Q_{\theta,k}(t)$ is continuously differentiable from
$\mathbb{R}$ to $\mathcal{L}(\mathcal{B}_{\gamma},\mathcal
{B}_{\gamma'})$, and for all $t\in\R$, $(\mathrm{d}Q_{\theta,k}/\mathrm{d}t)(t)$ is
the operator in $\mathcal{L}(\mathcal{B}_{\gamma},\mathcal
{B}_{\gamma'})$ associated to the kernel $Q_{\theta,k+1}(t)$.
\end{longlist}
Finally, we have ${\cal Q}_{k,\gamma,\gamma'} := \sup\{\|Q_{\theta
,k}(t)\|_{\gamma,\gamma'},\theta\in\Theta, t\in\R\}<\infty$,
and ${\cal Q}_{k,\gamma,\gamma'}$ depends on~$\xi$ but only via the
constant $C_\xi$ of  (\ref{d-m}).
\end{alem}

\begin{pf}
Set $\Delta_{\theta,k} := Q_{\theta,k}(t) - Q_{\theta,k}(t_0)$, and
let $0<\varepsilon\leq1$ be such that $\gamma+ (k+\varepsilon)/m
\leq\gamma'$. Using $|\mathrm{e}^{\mathrm{i}a}-1|\leq2|a|^{\varepsilon}$ and (\ref
{d-m}), we obtain for $f\in\cB_{\gamma}$:
\begin{eqnarray*}
|\Delta_{\theta,k}f(x)| &\leq& 2  |t-t_0|^{\varepsilon}  \|f\|
_{\gamma}  \int|\xi(\alpha_0,x,y)|^{k+\varepsilon}  V(y)^\gamma
Q_\theta(x,\mathrm{d}y) \\
&\leq& 2^{1+\fracb{k+\varepsilon}{m}}  {C_\xi}^{\fracb
{k+\varepsilon}{m}}  |t-t_0|^{\varepsilon}  \|f\|_{\gamma}
 \bigl(V^{\fracb{k+\varepsilon}{m}}(x)  Q_\theta V^\gamma(x) +
Q_\theta V^{\gamma'}(x) \bigr).
\end{eqnarray*}
Since the functions $V^{-\gamma}  Q_\theta V^\gamma$ and $V^{-\gamma
'}  Q_\theta V^{\gamma'}$ are bounded on $E$ uniformly in $\theta\in
\Theta$, we deduce that $\|\Delta_{\theta,k}f\|_{\gamma'} \leq
D_\xi  |t-t_0|^{\varepsilon}  \|f\|_{\gamma}$, where $D_\xi$ is a
positive constant\vadjust{\goodbreak} depending on~$C_\xi$ (but independent of $\theta$).
This gives (a). The proof of (b) is similar, using the operators
$Q_{\theta,k}(t) - Q_{\theta,k}(t_0) - (t-t_0)Q_{\theta,k+1}(t_0)$
and the inequality $|\mathrm{e}^{\mathrm{i}a}-1-\mathrm{i}a|\leq2|a|^{1+\varepsilon}$.
\end{pf}
\end{appendix}

%

\printhistory

\end{document}